\magnification=\magstep1
\hsize=16truecm
\input amstex
\TagsOnRight
\parindent=20pt
\parskip=1.5pt plus 1pt
\define\({\left(}
\define\){\right)}
\define\[{\left[}
\define\]{\right]}
\define\e{\varepsilon}

\define\supp {\sup\limits}

\define\summ{\sum\limits}
\define\prodd{\prod\limits}

\define\bigcupp{\bigcup\limits}

\centerline{\bf ESTIMATION OF WIENER--IT\^O INTEGRALS AND}
\centerline{\bf POLYNOMIALS OF INDEPENDENT GAUSSIAN}
\centerline{\bf RANDOM VARIABLES.}

\medskip
\centerline{\it P\'eter Major}

\centerline{\it Alfr\'ed R\'enyi Mathematical Institute}
\centerline{of the Hungarian Academy of Sciences}

\medskip
{\narrower{\narrower \noindent In this paper I prove good estimates 
on the moments and tail distribution of $k$-fold Wiener--It\^o 
integrals and also present their natural counterpart for polynomials 
of independent Gaussian random variables. The proof is based on the
so-called diagram formula for Wiener--It\^o integrals which yields
a good representation for their products as a sum of such integrals.
I intend to show in a subsequent paper that this method also yields
good estimates for degenerate $U$-statistics. The main result of this 
paper is a generalization of the estimates of
Hanson and Wright about bilinear forms of independent standard normal
random variables. On the other hand, it is a weaker estimate than the
main result of a paper of Lata{\l}a~[6]. But that paper contains an
error, and it is not clear whether its result is true. This question
is also discussed here.
\par}\par}

\medskip

\beginsection 1. Introduction. Formulation of the main results.

The goal of this paper is to give good estimates on the
tail-distribution and on high moments of Wiener--It\^o integrals.
This problem can be reformulated in an equivalent form to the
estimation of polynomials of independent Gaussian random variables.
The results obtained in such a way will be also proved.

This paper can be considered as a continuation of my investigation
in paper~[10] where good estimates were given for the
tail-distribution and high moments of Wiener--It\^o integrals and
(degenerate) $U$-statistics with the help of their variance. The
results of~[10] were proved by means of the so-called diagram
formula which yields a useful expression for the moments of
Wiener--It\^o integrals or degenerate  $U$-statistics. In the
present work it is shown that this method can be applied also in
cases when we have more information about a Wiener--It\^o integral
than its variance and we want to exploit this. 
I intend to prove similar improvements about the moments and 
tail-distribution of degenerate $U$-statistics
in a subsequent paper paper~[11].

Previous papers in this field (see [3] or [6]) dealt with the
estimation of polynomials of independent standard normal random
variables  and not of Wiener--It\^o integrals. I shall show that the
result of [3] about the estimation of bilinear forms of independent
standard normal random  variables is equivalent to the special case
of the main result in this  paper when only two-fold Wiener--It\^o
integrals are considered. On the other hand, paper~[6] formulates
sharper estimates about Gaussian polynomials of higher order than our
results. But the proof in [6] contains an error, hence some problems
arise with respect to this paper. I return to this question later.

First I introduce some notations needed to formulate the results of
the present paper.

Let us have a $\sigma$-finite non-atomic measure $\mu$ on some
measurable space $(X,\Cal X)$ together with a white noise $\mu_W$
with reference measure $\mu$, i.e. a set of jointly Gaussian random
variables $\mu_W(B)$ indexed by the sets $B\in\Cal X$ such that
$\mu(B)<\infty$, whose joint distribution is determined by the
relations $E\mu_W(B)=0$, $E\mu_W(A)\mu_W(B)=\mu(A\cap B)$ for all
sets $A, B\in\Cal X$ such that $\mu(A)<\infty$, $\mu(B)<\infty$. Let
us also introduce the quantity
$$
V_1^2(f)=\int f^2(x_1,\dots,x_k)\mu(\,dx_1)\dots \mu(\,dx_k) \tag1.1
$$
for a function $f$ of $k$ variables on the space $(X,\Cal X)$.
The $k$-fold Wiener--It\^o integral
$$
I_k(f)=\frac1{k!}\int f(x_1,\dots,x_k)\mu_W(\,dx_1)\dots\mu_W(\,dx_k)
\tag1.2
$$
can be defined for all functions $f$ such that $V_1(f)<\infty$.
(See e.g. [4] or~[8].) Here the knowledge of the definition of
Wiener--It\^o integrals is not assumed. We only need the so-called 
the diagram formula which enables us to calculate the moments of 
these random integrals. This result will be recalled in Section~3.

We are interested in good estimates on the probability
$P(k!|I_k(f)|>x)$ for large numbers $x>0$. This problem is closely
related to the question about good moment estimates $E(k!I_k(f))^{2M}$
for large values $M$. The results of paper [10] yield a good estimate
for these quantities with the help of $V_1(f)$. I recall them in the
following Theorem~A.

\medskip\noindent
{\bf Theorem A.} {\it Let a measurable space $(X,\Cal X)$ be given
together with a non-atomic $\sigma$-finite measure $\mu$ on it. Let $\mu_W$
be a white-noise with reference measure $\mu$, and let such a function
$f(x_1,\dots,x_k)$ of $k$ variables be given on the space $(X,\Cal X)$
for which $V_1^2(f)<\infty$ with the quantity $V_1(f)$ defined in~(1.1). 
Then the $k$-fold Wiener--It\^o integral  $I_k(f)$ introduced in~(1.2) 
satisfies the inequalities
$$
E(k!I_k(f))^{2M}\le C \(\frac {2kM}e\)^{kM}V_1(f)^{2M} \quad
\text{for all } M=1,2,\dots \tag1.3
$$
and
$$
P(k!|I_k(f)|>x)\le C\exp\left\{-\frac12\(\frac x{V_1(f)}\)^{2/k}\right\}
\quad \text{for all } x>0 \tag1.4
$$
with some appropriate universal $C>0$ depending only on the
multiplicity $k$ of the Wiener--It\^o integral.

This estimate is sharp in the following sense. There are
functions $f(x_1,\dots,x_k)$, (a function of the form
$f(x_1,\dots,x_k)=g(x_1)\dots g(x_k)$, $\int g^2(x)\mu(\,dx)<\infty$,
is an appropriate choice) for which the constant in the exponent of
the probability estimate (1.4) cannot be increased, i.e.\, the
inequality $P(k!|I_k(f)|>x)\le Ce^{-K(x/V_1(f))^{2/k}}$ does not
hold for $K>\frac12$.}
\medskip

Let me remark that $EI_k^2(f)=\frac1{k!}V^2_1(f)$ if $f$ is a function
symmetric in its variables (which may be assumed), and $EI_k(f)=0$.
So Theorem~A gave a good estimate on the moments and
tail-distribution of Wiener--It\^o integrals with the help of their
variance. If we have no more information about the kernel
function~$f$ of the Wiener--It\^o integral $I_k(f)$, than a bound on
$V_1(f)$ then we cannot improve the estimates in Theorem~A. On the
other hand, better estimates can be given with the help of some 
other appropriately defined quantities. In this paper such results 
will be proved. For this goal first I introduce some new quantities.

Given a finite set $K$, let $\Cal P=\Cal P(K)$ denote the set of all
partitions of this set $K$ to non-empty sets. For a finite set $K$ 
and a partition $P=\{A_1,\dots,A_s\}\in \Cal P(K)$ of this set let 
us define the class $\Cal F_P$ of appropriate sequences of functions
on the space $(X,\Cal X)$ by the following formula:
$$
\aligned
\Cal F_P&=\biggl\{g_r(x_j,\; j\in A_r),\; 1\le r\le s\colon\; \\
&\qquad \int g_r^2(x_j,\,j\in A_r)\prod_{j\in A_r}\mu(dx_j)\le1
\quad \text{for all } 1\le r\le s \biggr\}
\endaligned\tag1.5
$$
if $P=\{A_1,\dots,A_s\}\in\Cal P(K)$. This means that $\Cal F_P$
consists of a sequence of functions $g_r$ on $(X,\Cal X)$ whose
variables are indexed by the elements of corresponding sets $A_r$
in the partition~$\Cal P$, and with $L_2$-norm less than or equal
to 1 with respect to the appropriate product of the copies of the
measure~$\mu$.

Given a finite set $K$ and a function $f(x_j,\;j\in K)$ with 
arguments indexed by the elements of this set $K$ let us define for 
all partitions $P=\{A_1,\dots,A_s\}\in\Cal P= \Cal P(K)$ of the 
set $K$ the quantity
$$
V_P(f)=\sup_{(g_1,\dots,g_s)\in\Cal F_P}
\int f(x_j,\;j\in K)\prod_{1\le r\le s}g_r(x_j,\, j\in A_r)
\prodd_{j\in K}\mu(\,dx_j). \tag1.6
$$
Beside this, introduce the class $\Cal P_s=\Cal P_s(K)\subset\Cal P(K)$
of partitions of $K$ which consist of exactly $s$ elements, and put
$$
V_s(f)=\sup_{P\in\Cal P_s} V_{P}(f). \tag1.7
$$

I have defined the quantities $\Cal F_P$, $V_P(f)$ and $V_s(f)$
for a general finite set $K$ and square integrable function 
$f(x_j,j\in K)$, although in the formulation of our results they 
appear only in the special case $K=\{1,\dots,k\}$. But in the proofs 
we work with these notions in their general form. Let me also remark 
that the definition of $V_1(f)$ introduced in~(1.1) agrees with the 
definition of $V_s(f)$ in (1.7) with $s=1$ and $K=\{1,\dots,k\}$.  

The main result of this paper, an estimate about the moments and tail
distribution of a Wiener--It\^o integral $I_k(f)$, can be formulated
with the help of the quantities $V_s(f)$ introduced in formulas
(1.5),~(1.6) and~(1.7).

\medskip\noindent
{\bf Theorem about the tail-distribution of Wiener--It\^o
integrals.} {\it Let us consider a $k$-fold Wiener--It\^o integral
$I_k(f)$, $k\ge2$, defined in (1.1) by means of a white noise $\mu_W$
with a non-atomic $\sigma$-finite reference measure $\mu$ on a
measurable space $(X,\Cal X)$ and a measurable function
$f(x_1,\dots,x_k)$ of $k$ variables on the space $(X,\Cal X)$ such
that $V_1^2(f)<\infty$ for the quantity $V_1(f)$ defined in (1.1).
Then there exist some universal constants $C>0$, $C_1>0$ and $C_2>0$
depending only on the multiplicity $k$ of the Wiener--It\^o integral
$I_k(f)$ such that
$$
E(k!I_k(f))^{2M}\le C^MV_1(f)^{2M}\max\(M,
M^k\max_{2\le s\le k}\(\frac{V_s(f)}{V_1(f)}\)^{2/(s-1)}\)^M
\tag1.8
$$
for all $M=1,2,\dots,$, and
$$
\aligned
&P(|k!I_k(f)|>x)\\
&\qquad\le C_1 \exp\left\{-C_2 \min \( \frac{x^2}{V_1(f)^2},
\min_{2\le s\le k}
\(\frac x{V_1(f)^{1/(s-1)}V_s(f)^{(s-2)/(s-1)}}\)^{2/k}\)\right\}
\endaligned \tag1.9
$$
with the quantities $V_s(f)$, $1\le s\le k$, defined in formulas (1.5),
(1.6) and (1.7) with the set $K=\{1,\dots,k\}$.}

\medskip\noindent
{\it Remark.} The theorem about the tail-distribution of Wiener--It\^o
integrals may provide an essential improvement of Theorem~A in the case
when $V_s(f)$ is much smaller, than $V_1(f)$ for $2\le s\le k$. If we
have no information about the value of $V_s(f)$ for $s\ge2$, then we
can exploit the inequality $V_s(f)\le V_1(f)$ for $s\ge2$ and the
fact that inequalities (1.8) and (1.9) remain valid, if $V_s(f)$ on
its right-hand side is replaced by a larger number, for instance by
$V_1(f)$. In such a way we get slightly weaker estimates, than in 
Theorem~A. The estimates in both cases have the same structure, but
Theorem~A gives better information about the constants appearing
in them. Let me remark that although the estimate~(1.9) 
after the replacement of $V_s(f)$ by $V_1(f)$ has a form slightly 
different from the estimate~(1.4), this difference has no great
importance. It is not difficult to understand that in these 
estimates we may restrict our attention to the case $x\ge V_1(f)$, 
and in this case the term $\frac{x^2}{V_1(f)^2}$ can be dropped 
from the modified version of formula~(1.9).

\medskip
It will be more convenient to prove first the following simpler 
version of the above theorem and to deduce the result in the general 
case from it.

\medskip\noindent
{\bf Simplified version of the theorem about the tail distribution
of Wiener--It\^o integrals.} {\it Let a $k$-fold Wiener--It\^o
integral $I_k(f)$, $k\ge2$, with respect to a white noise $\mu_W$ 
with reference measure $\mu$ be given together with a real 
number~$R$, $0\le R\le 1$, in such a way that the kernel 
function~$f$ of the Wiener--It\^o integral and the number~$R$ 
satisfy the inequalities
$$
V_s(f)\le R^{s-1} \quad \text{for all } 1\le s\le k, \quad \text{and }
R\ge M^{-(k-1)/2} \tag1.10
$$
with some positive integer $M$. Then the inequality
$$
E(k!I_k(f))^{2M}\le C^M M^{kM}R^{2M} \tag1.11
$$
holds with this number~$M$ and some universal constant $C$ depending 
only on the multiplicity $k$ of the Wiener--It\^o integral.}

\medskip\noindent
{\it Reduction of the theorem about the tail-distribution of
Wiener--It\^o integrals to its simplified version.} Let us
introduce the function $\bar f=\frac{f}{V_1(f)}$, and put
$$
R=R_M=\max\( M^{-(k-1)/2}, 
\max_{1\le s\le k}\(\frac{V_s(f)}{V_1(f)}\)^{1/(s-1)}\).
$$
The function $\bar f$ and number $R$ satisfy the conditions
of the reduced theorem. Hence, this result implies that relation~(1.10)
holds with this function $\bar f$ and number $R$, which is equivalent
to relation~(1.8).

Formula~(1.9) can be proved in the standard way by means of 
formula~(1.8) and the Markov inequality
$P(|k!I_k(f)|>x)\le \frac{E(k!I_k(f))^{2M}}{x^{2M}}$
with a good choice of the parameter $M$. The choice of the closest
integer to
$$
\tilde C \min\(\frac{x^2}{V_1(f)^2},\min_{2\le s\le k}
\(\frac x{V_1(f)^{1/(s-1)}V_s(f)^{(s-2)/(s-1)}}\)^{2/k}\)
$$
for the parameter~$M$ with a sufficiently small $\tilde C>0$
(if $x>0$ is sufficiently large) supplies formula (1.9).

\medskip\noindent
{\it Remark.}\/ Formula (1.8) in the theorem about the 
tail-distribution of Wiener--It\^o integrals states in particular
that in the case $V_1(f)\le1$ the inequality $EI_k(f)^{2M}\le C^MM^M$
holds for all $M\le M_0^{2/(k-1)}$ with $M_0=\min\limits_{2\le s\le k}
V_s(f)^{-1/(s-1)}$, i.e. for such values $M$ the moments 
$EI_k(f)^{2M}$ have a bound similar to those of Gaussian random 
variables with expectation zero and variance smaller than a fixed 
positive number. This seems to be the most important part of this 
Theorem. Formula~(1.9) states a similar result about the tail 
distribution of $I_k(f)$. In the next session a result of 
Lata{\l}a~[6] will be discussed which states that the inequality 
$EI_k(f)^{2M}\le C^MM^M$ holds in a much larger interval, namely 
for $M\le M_0^2$. But the proof of this result contains an error, 
and it is not clear whether it holds.

The small value of the quantities $V_s(f)$, $2\le s\le k$, means
a sort of weak dependence property. This can be better understood
in the reformulation of our result for appropriate Gaussian 
polynomials of independent standard normal random variables, as 
it is done in the next section. In that reformulation some new 
quantities $\bar V_s(a(\cdot))$ defined with the help of the 
coefficients of the polynomials take the role of the numbers  
$V_s(f)$. The small value of these numbers $\bar V_s(a(\cdot))$ 
means that the random polynomial we consider is the sum of weakly 
dependent random variables. We have proved that even relatively
high moments of the some we consider behave like the moments of
Gaussian random variables under such conditions. It is an open 
question whether our result is sharp or higher moments of 
Wiener--it\^o integrals or random polynomials also satisfy a 
similar estimate  under the same conditions.

\medskip
This paper contains the proof of the above result. But before turning
to it I discuss what kind of estimates it yields for polynomials 
of independent standard normal random variables. This will be the
subject of Section~2. Beside this, this section contains a 
comparison of these estimates with earlier results in this field. 
In particular, a result of Lata{\l}a is discussed there together 
with the problem appearing in its proof. Section~3 contains a 
formula which helps to calculate the moments of a Wiener--It\^o 
integral. This formula is a consequence of the diagram formula for 
the product of Wiener--It\^o integrals, and it enables us to prove 
good moment estimates for Wiener--It\^o integrals if we can find 
good bounds on some integrals defined with the help of some 
diagrams. A closer study showed that it is useful to restrict our
attention to a special class of diagrams, to the so-called 
connected diagrams and to estimate the integrals related to them. 
A good bound on such integrals, called the Basic 
Estimate is also formulated in Section~3, and the 
proof of the main result of this paper is reduced to that of 
the Basic Estimate. In Section~4 a result called the Main 
Inequality is proved, and the Basic Estimate is proved with its 
help. Finally in Section~5 Lata{\l}a's result about an improvement 
of our estimates is discussed in more detail, and the question is 
investigated what kind of result has to be solved to decide 
whether it is true.

\beginsection 2. Bounds on random polynomials of Gaussian
random variables.

Let us take a natural counterpart of the estimation of $k$-fold 
Wiener--It\^o integrals, the estimation of some special polynomials 
of order~$k$ of independent standard normal random variables defined 
with the help of Hermite polynomials. In more detail, polynomials 
of the following form are considered. Let us have a sequence of 
independent, standard normal random variables $\xi_1,\xi_2\dots$, 
and introduce with their help the following random polynomials.
$$
Z_k=\sum\Sb((j_1,l_1),\dots,(j_s,l_s)),\; \\
j_u\neq j_{u'} \text{ if } u\neq u',\; l_1+\cdots+l_s=k \endSb
a((j_1,l_1),\dots,(j_u,l_u)) H_{l_1}(\xi_{j_1})\cdots H_{l_s}(\xi_{j_s}).
\tag2.1
$$
Here the coefficients $a((j_1,l_1),\dots,(j_u,l_u))$ are some real
numbers, and $H_l(x)$ denotes the Hermite polynomial of order~$l$
with leading coefficient~1. For the sake of simplicity let us assume
that the sum in formula (2.1) contains only finitely many terms. 
Infinite sums could also be allowed, but in that case some 
convergence problems should be handled.

It is more convenient to rewrite the random polynomials $Z_k$
in formula~(2.1) in a different form. This new version introduced
below means to work with Wick polynomials of Gaussian random
variables, i.e. to apply a multivariate generalization of
Hermite polynomials. To introduce the new representation of our
polynomials put
$$
H_{l_1}(\xi_{j_1})\cdots H_{l_s}(\xi_{j_s})=:\!
\undersetbrace l_1\text{-times} \to{\xi_{j_1},\dots,\xi_{j_1}},\dots,
\undersetbrace l_s\text{-times} \to{\xi_{j_s},\dots,\xi{j_s}} \!:.
$$
At the right-hand side of the last formula there is the product
of $k$ Gaussian random variables $\xi_{j_s}$, $1\le j_s\le n$,
$1\le s\le k$, between the two signs $:$, and some of the terms
$\xi_{j_s}$ may agree. For the sake of a more convenient notation
let us slightly extend the definition of the above expression. Let
us allow to write the terms in this product in an arbitrary order,
i.e.\ put
$$
:\! \xi_{j_{\pi(1)}},\dots,\xi_{j_{\pi(k)}}\!:
=:\! \xi_{j_1},\dots,\xi_{j_k}\!:
$$
for all permutations $\pi=(\pi(1),\dots,\pi(k))$ of the set
$\{1,\dots,k\}$. With such a notation the random polynomials we
are working with can be rewritten in the form
$$
Z_k=\sum a(n_1,\dots,n_k):\!\xi_{n_1}\xi_{n_2}\cdots\xi_{n_k}\!:
\tag2.2
$$
with some real coefficients $a(n_1,\dots,n_k)$ that can be calculated
by means of the coefficients $a((j_1,l_1),\dots,(j_u,l_u))$ in
formula~(2.1). In the subsequent considerations I shall estimate
random polynomials presented in the form~(2.2).

Let us consider the unit interval $[0,1]$ together with the Lebesgue
measure on it denoted by $\mu$, and let $\mu_W$ be a white noise with 
this reference measure~$\mu$. Take a complete orthonormal system of 
functions $\varphi_1(x),\varphi_2(x),\dots$ with respect to the 
Lebesgue measure on $[0,1]$ and put 
$\xi_n=\int\varphi_n(x)\mu_W(\,dx)$, $n=1,2,\dots$. Then 
$\xi_1,\xi_2,\dots$ is a sequence of independent, standard normal
random variables, and we do not change the distribution of the random
polynomial $Z_k$ in~(2.2) by choosing the above introduced standard
normal random variables $\xi_n$ in it. Thus it may be assumed that 
these standard normal random variables appear in the definition of 
the random polynomial~$Z_k$, and I shall exploit this liberty. With 
such a choice It\^o's formula for multiple Wiener--It\^o integrals 
(see e.g.~[4]) enables us to rewrite the random polynomial $Z_k$ in 
the form of a $k$-fold Wiener--It\^o integral. Such a representation 
of the random polynomial~$Z_k$ is useful, because it enables us to 
apply the results we know about Wiener--It\^o integrals in our 
investigations.

To find the above indicated representation observe that by It\^o's
formula
$$
:\!\xi_{n_1}\xi_{n_2}\cdots\xi_{n_k}\!:
=\int\varphi_{n_1}(x_1)\cdots\varphi_{n_k}(x_k)
\mu_W(\,dx_1)\dots\mu_W(\,dx_k),
$$
hence
$$
Z_k=k!I_k(f)=\int f(x_1,\dots,x_k)\mu_W(\,dx_1)\dots\mu_W(\,dx_k)
\tag2.3
$$
with
$$
f(x_1,\dots,x_k)=\sum a(n_1,\dots,n_k)
\varphi_{n_1}(x_1)\cdots\varphi_{n_k}(x_k).
\tag2.4
$$

We get an estimate for the moments and tail-distribution of
the random polynomial $Z_k$ with the help of the results formulated
in Section~1 if we can express the quantities $V_s(f)$, $1\le s\le k$,
for the function $f$ introduced in (2.4) by means of the coefficients
$a(n_1,\dots,n_k)$ in formula (2.2). We shall define some quantities
$\bar V_s$, $1\le s\le k$, with the help of the coefficients
$a(n_1,\dots,n_k)$ and show that they are equal to $V_s(f)$,
$1\le s\le k$, with the function $f$ defined in (2.4).

To define the desired quantities let us first introduce the set
$\Cal G_P$ corresponding to $\Cal F_P$ defined in~(1.5) for a partition
$P=\{A_1,\dots,A_s\}\in\Cal P(K)$ of the set $K=\{1,\dots,k\}$.
$$
\aligned
\Cal G_P&=\biggl\{b_r(n_j,\; j\in A_r),\; 1\le n_j<\infty
\text{ for all }j\in A_r,\; 1\le r\le s\colon\; \\
&\qquad \sum_{1\le n_j<\infty,\,1\le j\le r} b_r^2(n_j,\,j\in A_r)\le1
\quad \text{for all } 1\le r\le s \biggr\}
\endaligned
$$
if $P=\{A_1,\dots,A_s\}\in\Cal P(K)$. With the help of this notion
define similarly to the definition of $V_P$ in (1.6) the quantity
$$
\bar V_P(a(\cdot))=\sup_{(b_1,\dots,b_s)\in\Cal G_P}
\sum a(n_1,\dots, n_k)\prod_{1\le r\le s}b_r(n_j,\, j\in A_r) \tag2.5
$$
for a partition $P=\{A_1,\dots,A_s\}\in\Cal P(K)$ and a sequence
$a(n_1,\dots,n_k)$, $0\le n_j<\infty$ for all $1\le j\le k$ such
that $\sum a^2(n_1,\dots,n_k)<\infty$. Finally, put, similarly to
formula~(1.7),
$$
\bar V_s(a(\cdot))=\sup_{P\in\Cal P_s}\bar V_P(a(\cdot)). \tag2.6
$$

I claim that
$$
\bar V_P(a(\cdot))=V_P(f) \tag2.7
$$
for all partitions $P\in\Cal P(K)$, and as a consequence
$$
\bar V_s(a(\cdot))=V_s(f)   \tag2.8
$$
with the function $f$ defined in (2.4) and the sequence
$a(n_1,\dots,n_k)$ appearing in formula~(2.2).

To show that relation (2.7) holds fix some partition
$P=\{A_1,\dots,A_s\}\in\Cal P(K)$ with
$A_r=\{(j_1^{(r)},\dots,j_{l(r)}^{(r)})\}$, $1\le r\le s$, take 
a set of sequences
$$
(b_1(n_j,j\in A_1),\dots,b_s(n_j,j\in A_s))\in\Cal G_P,
$$ 
and correspond to it the set of functions $(g_1,\dots,g_r)$ defined
by the formula
$$
g_r(x_j,\; j\in A_r)= \!\!\!
\sum_{n_j^{(u)}\colon\; 1\le n_{j^{(r)}_u}<\infty, \; 1\le u\le l(r)} 
\!\!\! 
b_r(n_{j_1^{(r)}},\dots,n_{j_{l(r)}^{(r)}})
\prod_{u=1}^{l(r)}\varphi_{n_{j_u^{(r)}}}(x_{j_u^{(r)}}),
\quad 1\le r\le s. \tag2.9
$$
Then it follows from the Parseval formula that
$$
\int g^2_r(x_j,\;j\in A_r)\prod_{j\in A_r}\mu(\,dx_j)
=\sum_{n_j\colon\;1\le n_j<\infty,\; j\in A_r}  
b_r^2(n_j, j\in A_r), \quad 1\le r\le s,
$$
for the functions $(g_1,\dots,g_r)\in\Cal F_P$, and the mapping 
$b_r\to g_r$, $1\le r\le s$, defined in formula (2.9) is a 
one-to-one map from $\Cal G_P$ to $\Cal F_P$. Beside this, the 
Parseval formula also implies that
$$
\aligned
\sum a(n_1,\dots,n_k)\prod_{r=1}^s b_r(n_j,\;j\in A_r)
=\int f(x_1,\dots,x_k)\prod_{r=1}^s g_r(x_j,\,j\in A_r)
\prod_{j=1}^k\mu(\,dx_j)
\endaligned \tag2.10
$$
with the sequences $a(n_1,\dots,n_k)$ in (2.2) and the function
$f(x_1,\dots,x_k)$ in (2.4).

Then taking the supremum of both sides of formula (2.10) for all
$(b_1,\dots,b_r)\in\Cal G_P$, and exploiting that the two suprema
equal we get relation~(2.5). Indeed, the supremum of the left-hand
side equals $\bar V_P((a(\cdot)$ by definition. The supremum of
the right-hand side equals $V_P(f)$. Indeed, the properties of
the above defined one-to-one map imply that the supremum of
the right-hand side equals the supremum of the same expressions
if the supremum is taken for all $(g_1,\dots,g_s)\in\Cal F_P$. 
Finally, taking supremum in formula (2.5) for all $P\in\Cal P_s$ 
we get formula~(2.6).

The above considerations together with the theorem about the
tail-distribution of Wiener--It\^o integrals formulated in Section~1
yield the following result.

\medskip\noindent
{\bf Theorem about the estimation of moments and tail distribution
of polynomials of independent standard normal random variables.} {\it
Let us consider the random polynomial $Z_k$ of standard normal random
variables defined in formula~(2.2).
There exist some universal constants $C>0$, $C_1>0$ and $C_2>0$
depending only on the order $k$ of this random polynomial such that
$$
EZ_k^{2M}\le C^MV_1(a(\cdot))^{2M}\max\(M,M^k\max_{2\le s\le k}
\(\frac{\bar V_s(a(\cdot))}{\bar V_1(a(\cdot))}\)^{2/(s-1)}\)^M \tag2.11
$$
for all $M=1,2,\dots,$, and
$$
\aligned
P(|Z_k|>x)\le C_1 \exp
&\biggl\{-C_2\min\biggl(\frac{x^2}{\bar V_1(a(\cdot))^2}, \\
& \qquad \min_{2\le s\le k}\(\frac x{\bar V_1(a(\cdot))^{1/(s-1)}
\bar V_s(a(\cdot))^{(s-2)/(s-1)}}\)^{2/k}\biggr)\biggr\}
\endaligned \tag2.12
$$
with the quantities $\bar V_s(f)$, $1\le s\le k$, defined in formulas
(2.5) and (2.6).}

\medskip
The paper of Hanson and Wright~[3] contains some estimates on the tail
distribution of a bilinear form
$$
S_n=\sum_{i,j=1}^n a(i,j)(\xi_i\xi_j-E\xi_i\xi_j), \tag2.13
$$
where $\xi_1, \dots,\xi_n$ are independent standard normal random
variables, and $A=(a(i,j))$, $1\le i,j\le n$, is an $n\times n$
symmetric matrix. Hanson and Wright introduced the Hilbert--Schmidt
norm $\Lambda^2=\summ_{i,j=1}^n a(i,j)^2$ and the usual norm
$\|A\|=\supp_{|x|=1}|Ax|$ of the matrix $A=(a(i,j))$, where $|x|$
denotes the usual Euclidean norm of the vector $x=(x_1,\dots,x_n)$.
They proved the following estimate for the tail distribution of $S_n$.
$$
P(S_n>x)\le C_1\exp
\left\{-\min\(\frac{ C_2x}{\|A\|},\frac{C_2x^2}{\Lambda^2}\)\right\}
\tag2.14
$$
for all $x>0$ with some universal constants $C_1>0$ and $C_2>0$.

This inequality agrees with the estimate (2.12) in the theorem about
the estimation of moments and tail distribution of polynomials of
independent standard normal random variables in the case $k=2$.
Indeed, the random polynomials $Z_k$ defined in (2.2) agree
with the polynomials $S_n$ in (2.13), if they contain the
same coefficients $a(i,j)$. Beside this, in this case
$\bar V_1(a(\cdot))^2=\Lambda^2$,  and $\bar V_2(a(\cdot))=\|A\|$
for the terms
$\bar V_1(a(\cdot))$ and $\bar V_2(a(\cdot))$ defined in (2.5)
and~(2.6). The condition that the matrix $A=((a(i,j))$ must be
symmetric does not mean a real restriction, because a general
polynomial $S_n$ can be replaced by its symmetrization, which does
not change its distribution. With such a notation the estimate (2.14)
agrees with the estimate~(2.12) for $k=2$.

\medskip
Lata{\l}a (see~[6]) studied the estimation of polynomials of
independent standard normal random variables which have the
following special form.
$$
Z_k=\sum_{1\le j_s\le n, \,1\le s\le k} a(j_1,\dots,j_k)
\xi_{j_1}^{(1)}\dots\xi_{j_k}^{(k)},   \tag2.15
$$
where $\xi^{(s)}_1,\dots,\xi^{(s)}_n$, $1\le s\le k$, are
independent random sequences of independent standard normal random
variables. In this case he formulated an estimate, sharper than
ours. But the proof of his result contains an error, and it is
not clear whether it holds. Hence I present this inequality
as a conjecture. Its formulation presented here is
slightly different from that of paper~[6], but they are equivalent.
They have a similar relation to each other as the original and 
simplified versions of the theorem about the tail distribution of 
Wiener--It\^o integrals in Section~1.

\medskip\noindent
{\bf Lata{\l}a's conjecture.} {\it Let the coefficients of the
random Gaussian polynomial $Z_k$ of order $k$ defined in (2.15)
satisfy the inequality $\bar V_s(a(\cdot))\le M^{-(s-1)/2}$ with
some positive integer~$M$ for all $1\le s\le k$, with the
quantities $\bar V_s(a(\cdot))$ defined in formulas~(2.5)
and~(2.6). Then there exists some universal constant $C$
depending only on the order~$k$ of the random polynomial $Z_k$
such that
$$
EZ_k^{2M}\le C^M M^M. \tag2.16
$$
}

It follows from a result of de la Pe\~na and Montgomery--Smith~[1]
that if Lata{\l}a's conjecture holds for random polynomials of the
form~(2.14), then it also holds for polynomials of the form
$$
Z_k=\sum_{1\le j_s\le n, \,1\le s\le k,\; j_s\neq j_{s'}
\text{ if }s\neq s'}   a(j_1,\dots,j_k)\xi_{j_1}\dots\xi_{j_k},
$$
where $\xi_1,\dots,\xi_n$ are independent standard normal random
variables. With the help of this observation and some additional
work Lata{\l}a's conjecture can be verified for general polynomials
of the form~(2.2), provided that it holds in its original form. I
omit the details.

It is not difficult to see that Lata{\l}a's conjecture formulates
a sharper estimate than inequality~(2.11) in
the theorem about the estimation of moments and tail distribution
of polynomials of independent standard normal random variables. 
Indeed, relation~(2.16) implies that $Z_k^{2M}\le C^M M^M A^M$ for a
general polynomial of the form~(2.2) with
$A=\max\limits_{1\le s\le k}M^{s-1}\bar V_s^2(a(\cdot))$. To prove 
that this is a sharper inequality than formula~(2.11) it is enough 
to check that 
$$
A\le\max\(\bar V_1^2(a(\cdot)), 
\max\limits_{2\le s\le k}M^{k-1}\bar V_1^{2(s-2)/(s-1)}(a(\cdot))
V_s^{2/(s-1)}(a(\cdot))\)
$$ 
for this number~$A$. But this inequality clearly holds, because
$$
M^{s-1}\bar V_s^2(a(\cdot))\le M^{k-1}\bar V_1^{2(s-2)/(s-1)}a(\cdot))
V_s^{2/(s-1)}(a(\cdot))
$$ 
for all $2\le s\le k$, and the corresponding estimation for $s=1$ 
also holds.

\medskip
Lata{\l}a's argument heavily exploited the special form of the 
random polynomials in~(2.15). His method strongly exploited that
the terms in the sum~(2.15) are products 
$\xi_{j_1}^{(1)}\dots\xi_{j_k}^{(k)}$ with elements from 
independent copies of a random sequences. This made possible 
the application of a conditioning argument and the reduction of
the original problem to the estimation of the supremum of an
appropriately defined class of Gaussian random variables with its
help. But the estimation of such a supremum is very hard,
and at this point some serious difficulties arise. A
problem of the following type has to be considered.

Let us have a set of (jointly) Gaussian random variables, $\eta(x)$, 
$E\eta(x)=0$, $E\eta^2(x)\le1$, $x\in X$, indexed by a parameter 
set $X$, and try to give a good estimate on the expected value 
$E\(\supp_{x\in X}\eta(x)\)^{2M}$ for large positive integers~$M$. 
To study this problem let us introduce the (pseudo)metric 
$\rho_\alpha$ defined by the formula
$\rho_\alpha(x,y)=\(E(\eta(x)-\eta(y))^2\)^{1/2}$,
$x,y\in X$, in the parameter space~$X$. There is a natural way to
give good estimates on the moments of the supremum we are interested
in if we can give for all $\e>0$ a good estimate on the minimal number
$N(X,\rho_\alpha,\e)$ of balls of radius $\e$ with respect to the
distance $\rho_\alpha$ which cover the whole parameter set~$X$.

This number $N(X,\rho_\alpha,\e)$ can be estimated in the following
special case. If a probability measure $\mu_\e$ can be introduced
in the parameter set $X$ (or on an extension of the set~$X$, as it is 
done in Lata{\l}a's paper) for all $\e>0$ in such a way that 
$\mu_\e(y\colon\; \rho(x,y)>\e/2)\ge H(\e)$ with some function 
$H(\cdot)$ for all $x\in X$, then the inequality
$N(X,\rho_\alpha,\e)\le \frac1{H(\e)}$ holds. Lata{\l}a claimed that
such a construction is possible in the problem he is investigating.
He formulated two lemmas, Lemma~1 and Lemma~2 in his paper which
supply a good estimate, presented in Corollary~2, on
$N(X,\rho_\alpha,\e)$.

However, the proof of these Lemmas~1 and~2 is problematic. Lemma~1
contains a small inaccuracy (it states the upper bound $e^{-t^2/2}$
instead of $e^{-1/2t^2}$, and this wrong formula is written rather
consequently), but this seems to be a corrigible error. The main
problem is that Lemma~1 yields a too weak estimate which is not
sufficient to prove Lemma~2. In the explanation of this point I 
refer to the notation of paper~[6].

The right formula in the first line of the proof of Lemma~2 for $d=1$
(page~2319 in paper~[6]) would be
$$
\align
B_\alpha(\bold x, W_d^{[\bold x]}(\alpha,4t))
&=\{y\in R^{n_1}\colon\; \alpha(x-y)
\le W_d^{[\bold x]}(\alpha,4t)\}\\
&=\{y\in R^{n_1}\colon\; \alpha(x-y)\le4t E\alpha (xG_{n_1})\},
\endalign
$$
where $xG_{n_1}=(x_j g_j$, $1\le j\le n_1)$, with a standard normal
random vector $G_{n_1}=(g_1,\dots,g_{n_1})$. Hence in the proof of 
Lemma~2 for $d=1$ Lemma~1 should be applied
with the parameter $\bar t$ determined by the equation
$4\bar tE\alpha(G_{n_1})=4tE\alpha(xG_{n_1})$ (instead of the
parameter~$t$, as it is done in~[6]). This number $\bar t$ can be 
very small, since such vectors $x$ have to be considered for which
$\summ_{j=1}^{n_1}x_j^2\le1$. Hence Lemma~1 does not supply a
good estimate in such a case. In particular, the estimate we get,
depends on the dimension $n_1$ of the space $R^{n_1}$, i.e. of the
space where the random vector $G_{n_1}$ takes its values. On the other
hand, we need such estimates which do not depend on this dimension.

The question arises whether the proof can be saved despite of this
error. The hardest problem about Lata{\l}a's proof is hiding behind 
this question. In Section~5 I return to it. For the sake of simpler 
notations I shall consider only the case $k=3$. I show that 
Lata{\l}a's conjecture is equivalent to a rather hard estimate on 
the expected value of the supremum of some random multilinear forms 
whose study demands new ideas.

\medskip\noindent
{\bf 3. The diagram formula for Wiener--It\^o integrals.}

\medskip\noindent
In this section the diagram formula for products of Wiener--It\^o
integrals is formulated, and it is shown how the proof
of the simplified version of the theorem about the tail distribution
of Wiener--It\^o integrals can be reduced with its help to an
estimate that I call the Basic Estimate.

The diagram formula makes possible to rewrite the product of
Wiener--It\^o integrals as a sum of such integrals, and
as a consequence, it supplies a formula about the moments of
Wiener--It\^o integrals. It was shown in~[9] that this formula
yields a good estimate on the moments of Wiener--It\^o integrals.
In this paper it will be shown that if the quantities $V_s(f)$,
$1\le s\le k$, defined in~(1.6) and~(1.7) are very small, then
this method yields a better moment estimate. I recall this formula. 
It is the same result that I presented in paper~[9], only I 
made some small changes in the notation. The indices of the
arguments of the functions we are working with will be indexed in
a different way, because this simplifies the discussion.

The following problem is considered:
Let us have $m$ such real-valued functions
$$
f_j(x_{1},\dots,x_{k_j}), \quad 1\le j\le m,
$$
with $k_j$ variables on a measure space
$(X,\Cal X,\mu)$ with some $\sigma$-finite non-atomic measure
$\mu$ for which
$$
\int f_j^2(x_1,\dots,x_{k_j})
\mu(\,dx_1)\dots\mu(\,dx_{k_j})
<\infty \quad \text{ for all } 1\le j\le m \tag3.1
$$
together with a white noise $\mu_W$ with reference measure~$\mu$ 
on the space $(X,\Cal X)$. The Wiener--It\^o integrals
$k_j!I_{k_j}(f_j)$, $1\le j\le m$, introduced in (1.3) can be
defined with the above kernel functions $f_j$ and white noise
$\mu_W$. We are interested in a good explicit formula for the
expectation $E\(\prodd_{j=1}^m I_{k_j}(f_j)\)$. This formula,
which is a simple consequence of the diagram formula for
products of Wiener--It\^o integrals (see [2] or~[8]) will be 
presented below.

The expectation of the above product can be expressed by means of
some (closed) diagrams introduced below. A class of (closed)
diagrams denoted by $\bar\Gamma=\bar\Gamma(k_1,\dots,k_m)$ is defined
in the following way. A diagram $\gamma\in\bar\Gamma(k_1,\dots,k_m)$
consists of vertices of the form $(j,l)$, $1\le j\le m$,
$1\le l\le k_j$, and edges $((j,l),(j',l'))$, $1\le j,j'\le m$,
$1\le l\le k_j$, $1\le l'\le k_j'$. The set of vertices 
$(j,l)$, $1\le l\le k_j$, with a fixed number $j$ is called the 
$j$-th row of the diagram. All edges $((j,l),(j',l'))$ of a 
diagram $\gamma\in\bar\Gamma(k_1,\dots,k_k)$ connect vertices from
different rows, i.e. $j\neq j'$. It is also demanded that exactly
one edge starts from all vertices of a (closed) diagram $\gamma$.
The class $\bar\Gamma(k_1,\dots,k_m)$ of (closed) diagrams contains
the diagrams $\gamma$ with the above properties. Beside this, I
introduce the set $U(\gamma)$ containing the edges of the diagram
$\gamma$ for all $\gamma\in\bar\Gamma(k_1,\dots,k_m)$ and enumerate
their elements. An arbitrary enumeration is allowed, it
is only demanded that different edges must get different labels.
Let $N(\gamma)$ denote the set of indices of the edges in
$U(\gamma)$. In this section I shall consider only closed
diagrams in which every vertex is the end-point of some edge. In
the next section we have to work with more general, not necessarily
closed diagrams. I introduce their definition there.

Let us fix a diagram $\gamma\in\bar\Gamma(k_1,\dots,k_m)$. I
introduce the following function $u_\gamma(\cdot)$ on the vertices 
$(j,l)$, $1\le j\le m$, $1\le l\le k_j$.

For each vertex of $\gamma$ there is a unique
edge $((j,l),(j',l'))\in U(\gamma)$, i.e. an edge in
$\gamma$ which connects $(j,l)$ with some other vertex $(j',l')$ of
the diagram. If this edge has label $n\in N(\gamma)$, then we put
$u_\gamma((j,l))=n$.

Given a fixed diagram $\gamma\in\bar\Gamma(k_1,\dots,k_m)$ let us
rewrite the functions $f_j$ with reindexed variables as
$f_j(x_{u_\gamma(j,1)},\dots,x_{u_\gamma(j,k_j)})$, $1\le j\le m$,
with the help of the above defined function $u_\gamma(\cdot)$.
(Two variables get the same index if the vertices related to them
were connected by an edge of the diagram $\gamma$.) Define the
product of these reindexed variables
$$
\bar F_\gamma(x_n,\; n\in N(\gamma))=\prod_{j=1}^m
f_j(x_{u_\gamma(j,1)},\dots,x_{u_\gamma(j,k_j)}) \tag3.2
$$
together with the integral of these functions
$$
F_\gamma=F_\gamma(f_1,\dots,f_m)=\int
\bar F_\gamma(x_n,\; n\in N(\gamma))\prod_{n\in N(\gamma)}\mu(\,dx_n)
\tag3.3
$$
for all $\gamma\in \bar\Gamma(k_1,\dots,k_m)$.

The expected value of the product of Wiener--It\^o integrals
$k_j!I_{k_j}(f_j)$, $1\le j\le m$, can be expressed with the
help of the above quantities $F_\gamma$ in the following way.

\medskip\noindent
{\bf Formula about the expected value of products of Wiener--It\^o
integrals.} {\it Let us consider the Wiener--It\^o integrals
$k_j!I_{k_j}(f_j)$ of some functions $f_j$, $1\le j\le m$,
satisfying relation (3.1) with respect to a white noise $\mu_W$
with reference measure $\mu$. The expected value of this product
satisfies the identity
$$
E\left(\prod_{j=1}^m k_j!I_{k_j}(f_j)\right)
=\sum_{\gamma\in\bar\Gamma(k_1,\dots,k_m)} F_\gamma \tag3.4
$$
with the numbers $F_\gamma$ defined in (3.2) and (3.3).}
\medskip

To get a good estimate on the expectation of the product of
Wiener--It\^o integrals by means of formula (3.4) we need a
good bound on the quantities $F_\gamma$. For this goal it is
useful to rewrite them by means of an appropriate recursive 
formula. To present such a formula let us define the
restrictions $\gamma^r$, $1\le r\le m$, of all (closed)
diagrams $\gamma\in\bar\Gamma(k_1,\dots,k_m)$ to its first $r$
rows, $1\le r\le m$.

More explicitly, the diagram $\gamma^r$ contains the vertices $(j,l)$,
$1\le j\le r$, $1\le l\le k_j$ and those edges of $\gamma$ whose
end-points are vertices in one of the first $r$ rows in $\gamma$.
I shall call $\gamma^r$ also a diagram, although it may have
vertices from which no edge starts. I return to this point in
the next section. Let $U_1(\gamma^r)$ denote the set of edges
and $U_2(\gamma^r)$ the set of those vertices of $\gamma^r$ 
from which no edge starts in $\gamma^r$, i.e. which are 
connected with a vertex $(j',l')$ with $j'>r$ in~$\gamma$. 
Let all vertices of $\gamma^r$ get the same enumeration they 
got as a vertex of $\gamma$. Let $N_1(\gamma^r)$ denote the 
set of indices of the vertices which are end-points of an edge 
in $U_1(\gamma^r)$, and $N_2(\gamma^r)$ the set of indices of 
the vertices in $U_2(\gamma^r)$. Put 
$N(\gamma^r)=N_1(\gamma^r)\cup N_2(\gamma^r)$.

Let us define, similarly to the quantities $\bar F_\gamma$ and
$F_\gamma$, the functions
$$
\bar F_{\gamma^r}(x_n,\; n\in N(\gamma^r))=\prod_{j=1}^r
f_j(x_{u_\gamma(j,1)},\dots,x_{u_\gamma(j,k_j)}) \tag3.5
$$
and
$$
F_{\gamma^r}(x_n,\;n\in N_2(\gamma^r))=\int
\bar F_{\gamma^r}(x_n,\;n\in N(\gamma^r))
\prod_{n\in N_1(\gamma_r)}\mu(\,dx_n)
\tag3.6
$$
for all $1\le r\le m$, and $\gamma\in\bar\Gamma(k_1,\dots,k_m)$. In
the case $N_1(\gamma^r)=\emptyset$ the integral at the right-hand
side equals $\bar F_{\gamma^r}(x_n,\;n\in N(\gamma^r))$. In general,
the following convention is applied. If we integrate a function
with respect to a product measure, then in the case when this
product contains zero terms, then the integral equals the
function itself.

It is not difficult to check that
$$
F_{\gamma^1}(x_n,\;n\in N_2(\gamma^1))
=f_1(x_{u_\gamma(1,1)},\dots,x_{u_\gamma(1,k_1)}),
\tag3.7
$$
and
$$
\aligned
&F_{\gamma^r}(x_n,\;n\in N_2(\gamma^r))\\
&\qquad =\int F_{\gamma^{r-1}}(x_n,\;n\in N_2(\gamma^{r-1}))
f_r(x_{u_\gamma(r,1)},\dots,x_{u_\gamma(r,k_r)}) \!\!\!
\prod_{n\in N_2(\gamma^r)\setminus N_2(\gamma^{r-1})}
\!\!\! \mu(\,dx_n)
\endaligned \tag3.8
$$
for all $2\le r\le m$. Beside this,
$$
F_\gamma=F_{\gamma^m}(x_n,\;n\in N_2(\gamma^m)). \tag3.9
$$
(Actually, $N_2(\gamma^m)=\emptyset$.)
Relations (3.7), (3.8) and (3.9) yield a recursive formula for
the quantity $F_\gamma$.

Next I introduce the notion of connected (closed) diagrams. They
turned out to be a useful object, because the quantity $F_\gamma$
can be better bounded for a connected than for a general
diagram $\gamma$, and it is enough to bound them to prove our
results.

\medskip\noindent
{\bf Definition of connected (closed) diagrams.} {\it A (closed)
diagram $\gamma\in\bar\Gamma(k_1,\dots,k_m)$ is connected if for all
sets of rows $A\subset\{1,\dots,m\}$ of the diagram $\gamma$ such
that $1\le|A|\le m-1$ there is such an edge
$((j_1,l_1),(j_2,l_2))\in U(\gamma)$ of the diagram $\gamma$ for
which $j_1\in A$ and $j_2\notin A$.}
\medskip

The following result, called the Basic Estimate yields a bound on
$F_\gamma$ for connected diagrams~$\gamma$.
This estimate turned out to be sufficient for our purposes.

\medskip\noindent
{\bf Basic Estimate.} {\it Let us consider a\/ {\rm connected},
closed diagram $\gamma\in\bar\Gamma(k_1,\dots,k_m)$, $m\ge2$,
and some functions $f_j$ of $k_j$ variables on a measure space
$(X,\Cal X,\mu)$, $1\le j\le m$, which satisfy the inequality
$V_s(f_j)\le R^{s-1}$ with some $0\le R\le1$, for all
$1\le j\le m$ and $1\le s\le k_j$. (The quantities $V_s(f)$ were
introduced in formulas (1.6) and (1.7).) The quantity
$F_\gamma=F_\gamma(f_1,\dots,f_m)$, introduced in~(3.3) satisfies
the inequality
$$
|F_\gamma|=|F_\gamma(f_1,\dots,f_m)|\le R^{m-2}. \tag3.10
$$
}

I show how the {\it simplified version of the theorem about the
tail distribution of Wiener--It\^o integrals}\/ can be proved with
the help of the Basic Lemma. To do this first I show that all
diagrams $\gamma\in\bar\Gamma(k_1,\dots,k_m)$ can be decomposed to
the union of disjoint connected diagrams in a unique way, and the
quantity $F_\gamma$ equals the product of the numbers
$F_{\{\cdot\}}$ corresponding to these connected diagrams.

More explicitly, there is a unique partition $A_1=A_1(\gamma)$,\dots,
$A_u=A_u(\gamma)$ of the set of rows $\{1,\dots,m\}$ of
$\gamma\in\bar\Gamma(k_1,\dots,k_m)$ in such a way that the diagram
$\gamma$ equals the union of the diagrams $\gamma_{A_r}$,
$1\le r\le u$, where $\gamma_{A_r}$ is the restriction of $\gamma$
to the rows in $A_r$, i.e. it contains the rows of $\gamma$ with
indices in the set $A_r$ together with the edges connecting
vertices from these rows. Beside this, all diagrams $\gamma_{A_r}$
of this decomposition must be connected. Also the identity
$$
F_\gamma=\prod_{r=1}^u F_{\gamma_{A_r}} \tag3.11
$$
holds.

In more detail, the restriction $\gamma_{A_r}$ of $\gamma$ to the
rows with indices in $A_r$ consists of vertices in the rows with
index in the set of $A_r$, and those edges whose both end-points
are among these vertices. The decomposition of $\gamma$ to the
diagrams $\gamma_{A_r}$, $1\le r\le u$, also means that $\gamma$
has no such edge which connects vertices from $\gamma_{A_r}$ with
vertices from $\gamma_{A_{r'}}$ with some $r\neq r'$. The
connectedness of $\gamma_{A_r}$ means that for all subsets
$B\subset A_r$, $B\neq\emptyset$ and $B\neq A_r$ there is an edge
of $\gamma_{A_r}$ which connects a vertex in a row with index
in $B$ and a vertex in a row with index in $A_r\setminus B$.
The quantities $F_{\gamma_{A_r}}$ can be defined similarly to
$F_\gamma$, for instance by a natural adaptation of the recursive
formulas (3.7), (3.8) and~(3.9) in the calculation of $F_\gamma$
to the calculation of $F_{\gamma_{A_r}}$. In this adaptation we
can write similar recursive relations, only the row indices
$1,\dots,m$ must be replaced by rows with indices 
$v_1,\dots,v_{|A_r|}$ if $A_r=\{v_1,\dots,v_{|A_r|}\}$ with 
$v_1<v_2<\cdots<v_{|A_r|}$.
It is not difficult to check that the Basic Estimate also implies
that under the conditions of this result the inequality
$|F_{\gamma_A}|\le R^{|A|-2}$ also holds for any connected
(closed) diagram with rows in the set $A\subset\{1,\dots,m\}$.

To find the desired decomposition of a diagram
$\gamma\in\bar\Gamma(k_1,\dots,k_m)$ let us define the graph $G(\gamma)$
with vertices $1,\dots,m$ in which two vertices $j_1$ and $j_2$ are
connected with an edge if and only if the diagram $\gamma$ contains
an edge which connects two vertices from the $j_1$-th and $j_2$-th rows.
Let $A_1,\dots,A_u$ be the connected disjoint components of this graph.
Then it is not difficult to see that $\gamma_{A_1},\dots,\gamma_{A_u}$
supplies the desired decomposition of the diagram $\gamma$, and also
relation (3.11) holds.

Given a set $A\subset \{1,\dots,m\}$ let
$\bar{\Cal\Gamma}_{\text{c}}(k_j,\,j\in A)$ denote the class of
those connected (closed) diagrams whose rows are the sequences
$\{(j,1),\dots,(j,k_j)\}$, $j\in A$. Put
$$
K_{\text{c}}(A)=\sum_{\gamma_A\in
\bar{\Cal\Gamma}_{\text{c}}(k_j,\,j\in A)}F_{\gamma_A},
\tag3.12
$$
and let us also introduce for a partition $P=\{A_1,\dots,A_u\}\in\Cal
P(\{1,\dots,m\})$ of the set $\{1,\dots,m\}$ the class of diagrams
$\bar\Gamma(A_1,\dots,A_u)=\bar\Gamma(A_1,\dots,A_u|k_1,\dots,k_m)$
containing those diagrams  $\gamma\in\bar\Gamma(k_1,\dots,k_m)$ whose
decomposition to connected components consists of diagrams with
rows indexed by the sets $A_1$,\dots, and $A_u$, i.e. of diagrams
$\gamma_{A_r}\in\bar{ \Cal \Gamma}_{\text{c}}(k_j,\,j\in A_r)$,
$1\le r\le u$. Define for all partitions $\{A_1,\dots,A_u\}$ of the 
set $\{1,\dots,m\}$ the quantity
$$
K(A_1,\dots,A_u)=\sum_{\gamma\in\bar\Gamma(A_1,\dots,A_u)} F_\gamma.
$$

It is not difficult to see with the help of relations (3.11)
and~(3.12) that
$$
K(A_1,\dots,A_u)=K_{\text{c}}(A_1)\cdots K_{\text{c}}(A_u)
$$
for all partitions $\{A_1,\dots,A_u\}$ of the set $\{1,\dots,m\}$.
Summing up this identity for all classes $\bar\Gamma(A_1,\dots,A_u)$
we get with the help of the identity (3.4) and the fact that each
diagram has a unique decomposition to connected diagrams that
$$
E\left(\prod_{j=1}^m k_j!I_{k_j}(f_j)\right)
=\sum_{\{A_1,\dots,A_u\}\in\Cal P(\{1,\dots,m\})}
K_{\text{c}}(A_1)\cdots K_{\text{c}}(A_u) \tag3.13
$$
In the next calculations I shall restrict my attention to the case
$m=2M$, $f_j=f$, $k_j=k$ for all $1\le j\le 2M$. I give a good upper
bound on $EI_k(f)^{2M}$ with the help of relations (3.12), (3.13)
and the Basic Estimate for such parameters $R$ which satisfy the
relations $V_s(f)\le R^{s-1}$ for all $1\le s\le k$, and 
$R\ge M^{-(k-1)/2}$. In these calculations I also
exploit that a connected closed diagram has at least two rows,
hence only such partitions $\{A_1,\dots,A_u\}\in\Cal P(\{1,\dots,m\})$
have to be considered in the sum at the right-hand side of (3.13)
for which $|A_r|\ge2$ for all $1\le r\le u$.

First I show that if $V_s(f)\le R^{s-1}$ for all
$1\le s\le k$, then
$$
|K_{\text{c}}(A)|\le (k|A|)^{k|A|/2}R^{|A|-2} \quad\text{for all }
A\subset\{1,\dots,m\} \tag3.14
$$
for the quantity $K_{\text{c}}(A)$ defined in~(3.12). Indeed, the
Basic Estimate implies that $|F(\gamma_A)|\le R^{|A|-2}$ for all
terms in the sum at the right-hand side of~(3.12), and clearly there
are less than $(k|A|)^{k|A|/2}$ diagrams with $|A|$ rows and in each
rows $k$ vertices. (At this point we bounded the number of all (and
not only the number of all connected) diagrams.)

To estimate the expression at the right-hand side of (3.13)
let us introduce the class of those partitions
$\Cal P_{u,t_1,\dots,t_u}$, $t_1+\cdots+t_u=2M$, $t_r\ge2$,
$1\le r\le u$, of the set $\{1,\dots,2M\}$ which consist of $u$
sets $A_1,\dots,A_u$, and the set $A_r$ has $t_r$ elements,
$1\le r\le u$. (The set $\Cal F_{u,t_1,\dots,t_u}$ depends on the
number $u$ and the set $\{t_1,\dots,t_u\}$, but it does not depend on
the order of the  elements $t_r$ in this set.) Let us first estimate
the contribution of the partitions $\Cal P_{s,t_1,\dots,t_u}$ to the
sum in~(3.13). I claim that
$$
\aligned
&\left|\sum_{\{A_1,\dots,A_u\}\in\Cal P_{u,t_1,\dots,t_u}}
K_{\text{c}}(A_1)\cdots K_{\text{c}}(A_u)\right| \\
&\qquad\le |\Cal P_{u,t_1,\dots,t_u}|
\prod_{r=1}^u (kt_r)^{kt_r/2}R^{t_r-2}
\le C^M M^M\[M^{(k-1)}R^2\]^{M-u}
\endaligned\tag3.15
$$
with a universal constant $C>0$ depending only on the parameter~$k$.

The first inequality in (3.15) is a straight consequence of
relation~(3.14). To prove the second inequality we need a good
bound on $|\Cal P_{u,t_1,\dots,t_u}|$. To get it let us first list 
the element $A_1,\dots,A_u$ of a partition in
$\Cal P_{u,t_1,\dots,t_u}$ in the following way. Let $A_1$ be the
set which contains the number~1, $A_2$ the set containing the 
smallest number not contained in $A_1$, e.t.c..  Let $t_r$ denote 
the cardinality of the set $A_r$ with such an indexation. Then the 
number $|\Cal P_{u,t_1,\dots,t_u}|$ can be bounded in the following
way.
$$
|\Cal P_{u,t_1,\dots,t_u}|\le \prod_{r=1}^u \frac{(2M)^{t_r-1}}{(t_r-1)!}
\le C_1^M\frac{M^{2M-u}}{\prodd_{r=1}^u t_r^{t_r}} \tag3.16
$$
with some appropriate $C_1>0$.

The first inequality in relation (3.16) can be simply checked. To
prove the second inequality let us first observe that
$\prodd_{r=1}^ut_r\le \(\frac1u\summ_{r=1}^u t_r\)^u=\(\frac{2M}u\)^u
\le C_2^M$ with some universal constant $C_2$. This inequality together
with the Stirling formula imply that
$\prodd_{r=1}^u(t_r-1)!\ge C_2^{-M}\prodd_{r=1}^ut_r!
\ge C_3^{-M}\prodd_{r=1}^u t_r^{t^r}$ with some appropriate value
$C_3$, hence relation (3.16) holds.

It follows from relation (3.16) that
$$
|\Cal P_{u,t_1,\dots,t_u}|\prod_{r=1}^u (kt_r)^{kt_r/2}R^{t_r-2}
\le C_4^MM^{2M-u}\(\prodd_{r=1}^ut_r^{(k-2)t_r/2}\)R^{2M-2u}.
\tag3.17
$$
Let us consider the maximum of the right-hand side in the
inequality~(3.17) in the parameters $t_1,\dots,t_u$ with a fixed
value~$u$. Since $k-2\ge0$, this expression takes its maximum if
$t_1=2M-2u+2$ and $t_r=2$, $2\le r\le u$. Hence
$$
\align
&|\Cal P_{u,t_1,\dots,t_u}|\prod_{r=1}^u (kt_r)^{kt_r/2}R^{t_r-2}\\
&\qquad\le C_4^MM^{2M-u}(2M-2u+2)^{(k-2)(M-u+1)}
2^{(u-1)(k-1)}R^{2M-2u} \\
&\qquad \le C^MM^{2M-u}M^{(k-2)(M-u)}R^{2M-2u}
=C^MM^M\[M^{k-1}R^2\]^{M-u}
\endalign
$$
with some $C>0$ depending only on the parameter~$k$. This inequality
finishes the proof of relation~(3.15).

If the parameter $R$ satisfies also the inequality $R\ge M^{-(k-1)/2}$,
i.e. $M^{k-1}R^2\ge1$, then relation~(3.15) implies that
$$
\left|\sum_{\{A_1,\dots,A_s\}\in\Cal P_{s,t_1,\dots,t_s}}
K_{\text{c}}(A_1)\cdots K_{\text{c}}(A_s)\right|\le C^M M^M(M^{k-1}R^2)^M.
 \tag3.18
$$

Finally, I show that relations (3.13) (with $m=2M$ and $f_j=f$ for 
all $1\le j\le 2M$) and~(3.18) imply that
$$
E( k!I_{k}(f))^{2M}\le2^{2M}C^MM^{kM}R^{2M}
$$
if $V_s(f)\le R^{(s-1)/2}$ for all $1\le s\le k$, and $R\ge
M^{-(k-1)/2}$. This means that the simplified version of the theorem 
about the tail distribution of Wiener--It\^o integrals follows from 
the Basic Estimate.

To see that relations (3.13) and (3.18) really imply the last 
inequality it is enough to observe that the class of partitions 
$\Cal P(\{1,\dots,2M\})$ is the union of the classes of partitions 
$\Cal P_{u,t_1,\dots,t_u}$, and there are less than $2^{2M}$ classes 
$\Cal P_{u,t_1,\dots,t_u}$. (The number of such classes is bounded 
by the number of sets of positive integers $\{t_1,\dots,t_u\}$ such 
that $t_1+\cdots+t_u=2M$.)

\medskip
We could prove our result by means of a good estimate on the
quantity $F_\gamma$ for connected diagrams. The introduction of
connected diagrams was very important in these considerations,
because the Basic Estimate holds only for such diagrams. I show an
example of non-connected (closed) diagrams which satisfies only a
much weaker estimate. I shall consider an appropriate diagram
$\gamma\in\Gamma(k_1,\dots,k_m)$ and a function $f_j=f$ with $m=2M$,
$k_j=k$ for all $1\le j\le2M$. Such a function $f$ will be taken
which is symmetric in its variables, and $V_1(f)=1$. Let the diagram 
$\gamma$ I consider have the property that its rows can be put into 
pairs in such a way that edges can connect only vertices from rows 
which are paired. For such a diagram $F_\gamma=V_1(f)^{2M}=1$, and 
this value is much larger than the bound in the Basic Estimate.
But there are relatively few such diagrams, their number equals
$\frac{(2M)!}{2^MM!}(k!)^M$. Hence the relatively great value
of $F_\gamma$ for such diagrams causes no problem.

The estimation of the moments with the help of connected diagrams
corresponds to a classical method of probability theory, to the
estimation of moments by means of semi-invariants. The quantity
$K_{\text{c}}(A)$ introduced in formula~(3.12) is actually the
semi-invariant of the random variables $k_j!I_{k_j}(f)$ for $j\in A$,
and the identity (3.13) is a special case of the formula about
the expression of the expectation of products of random variables
by means of semi-invariants. The semi-invariants are estimated in
this paper with the help of the Basic Estimate which will be proved 
by means of a result called the Main Inequality in the next Section. 
In the application of this inequality we strongly exploit that we are
working with connected diagrams. Our approach shows some similarity
with the High Temperature Expansion in Statistical Physics. I do
not need the precise meaning of the notions and methods mentioned
in this paragraph, hence I omit their detailed discussion.
Some useful information about semi-invariants can be found in the 
second chapter of the book~[12].

\beginsection 4. The proof of the Basic Estimate.

The Basic Estimate will be proved by means of an inductive
proposition about the behaviour of the functions $F_{\gamma^r}$,
$1\le r\le m$, where $\gamma^r$ is the restriction of the closed
diagram $\gamma$ to its first $r$ rows. To formulate and prove this
result the notion of diagrams, introduced in Section~3, will be 
generalized. Such diagrams will be defined which may be not closed, 
i.e. which may have such vertices from which no edge starts. Some
other objects related to these new diagrams will be also introduced
and some results needed in our further discussion will be
formulated.

A diagram with rows indexed by a finite set $A=\{j_1,\dots,j_m\}$,
$1\le j_1<j_2<\cdots<j_m$ of the positive integers
and with row length $j_t$, $1\le t\le m$, is a graph whose
vertices are the pairs $(j_t,l)$ with $1\le t\le m$ and
$1\le l\le k_{j_t}$. The set of points
$\{(j_t,l)\colon\; 1\le l\le k_{j_t}\}$ is called the $t$-th row of
the diagram. A diagram may contain such edges $((j_t,l),(j_{t'},l'))$
for which $j_t,j_{t'}\in A$, $1\le l\le k_{j_t}$,
$1\le l'\le k_{j_{t'}}$, and $j_t\neq j_{t'}$. In words this means
that edges can connect only vertices from different rows. Beside this,
it is required that from each vertex there starts either zero or one
edge. Graphs satisfying all these properties will be called diagrams.
The class of diagrams with rows indexed by a set $A$ and with $k_j$
element in the row indexed by $j\in A$ will be denoted by
$\Gamma(k_j|j\in A)$. The main difference between diagrams and closed
diagrams introduced in the previous section is that a general diagram
may contain vertices from which no edge starts. Those vertices of a
diagram $\gamma$ from which no edge starts will be called open
vertices. The notion of connected diagrams will be also introduced
for this more general class of diagrams.

\medskip\noindent
{\bf Definition of connected diagrams in the general case.} {\it
A diagram $\gamma\in\Gamma(k_j|j\in A)$ is connected if
for all sets $B\subset A$ such that $B\neq\emptyset$ and $B\neq A$
there exists an edge $((j_t,l),(j_{t'},l'))$ of $\gamma$ such that
$j_t\in B$ and $j_{t'}\in A\setminus B$.}

\medskip
Similarly to closed diagrams, general diagrams
$\gamma\in\Gamma(k_j|j\in A)$ also have a unique decomposition to
connected diagrams. To formulate this statement precisely let us
first introduce the reduction of a diagram to some of its rows.
Given a diagram $\gamma\in\Gamma(k_j|j\in A)$, $A=\{j_1,\dots,j_m\}$
and a set $B\subset A$ let $\gamma_B\in\Gamma(k_j|j\in B)$ denote
that diagram whose vertices are points of the form $(j_t,l)$,
$j_t\in B$, $1\le l\le k_{j_t}$, and two vertices in $\gamma_B$ are
connected by an edge if and only if they are connected by an edge
in $\gamma$. With this notation we can state that for all diagrams
$\gamma\in\Gamma(k_j|j\in A)$ there exists a unique partition
$A_1,\dots,A_u$ of the set of rows $A$ in such a way that all
diagrams $\gamma_{A_r}$, $1\le r\le u$, are connected, and the
diagram $\gamma_A$ is the union of the diagrams $\gamma_{A_r}$,
$1\le r\le u$. This can be proved similarly to the case of closed
diagrams. One has to take the graph whose vertices are the points of
the set $A$ and draw an edge between two vertices $j_t\in A$,
and $j_{t'}\in A$ if there is an edge in the diagram $\gamma$
connecting some vertices from the $j_t$-th and the $j_{t'}$-th row.
By taking the decomposition of this graph to connected components we
also get the  decomposition of the diagram $\gamma$ to connected
components.

Given a diagram $\gamma\in\Gamma(k_j|j\in A)$ with
$A=\{j_1,\dots,j_m\}$ together with some functions
$f_t(x_{(j_t,1)},\dots,x_{(j_t,k_{j_t})})$, $1\le t\le m$, a
natural generalization of the quantities $\bar F_\gamma$ and
$F_\gamma$ defined in (3.2) and (3.3) can be introduced, and
it can be shown that they have similar properties.
To introduce these quantities let us enumerate first the edges 
and then the vertices of the diagram $\gamma$. A vertex from which 
an edge starts gets the same label as the edge starting from it. The 
remaining vertices, from which no edge starts get a new label. In 
this enumeration two vertices get the same index if and only if they 
are connected by an edge. Let $u_\gamma{(j,l)}$ denote the label of 
 the vertex $(j,l)$ in this enumeration of the 
vertices of a diagram $\gamma$. Let $N(\gamma)$ denote the set of 
labels of all vertices and $N_1(\gamma)$ the set of labels of all 
open vertices in $\gamma$. Then we define, similarly to 
formulas (3.2) and~(3.3) the quantities
$$
\bar F_\gamma(x_n,\; n\in N(\gamma))=\prod_{t=1}^m
f_t(x_{u_\gamma(j_t,1)},\dots,x_{u_\gamma(j_t,k_{j_t})}) \tag4.1
$$
and
$$
\aligned
F_\gamma(x_n,\; n\in N_1(\gamma))
&=F_\gamma(x_n,\;n\in N_1(\gamma)|f_1,\dots,f_m)\\
&=\int \bar F_\gamma(x_n,\; n\in N(\gamma))
\prod_{n\in N_2(\gamma)}\mu(\,dx_n)
\endaligned \tag4.2
$$
with $N_2(\gamma)=N(\gamma)\setminus N_1(\gamma)$.

Given a diagram $\gamma\in\Gamma(k_j|j\in A)$ with
$A=\{j_1,\dots,j_m\}$ let $\gamma^r$ denote its restriction to its
first $r$ rows, i.e. to rows in $A_r=\{j_1,\dots,j_r\}$, $1\le r\le m$.
The natural modification of the recursive relations~(3.7), (3.8)
and~(3.9) remain valid also for general diagrams~$\gamma$. The only
difference we have to make is to rewrite the indices $u_\gamma(1,l)$
and $u_\gamma(r,l)$ by $u_\gamma(j_1,l)$ and $u_\gamma(j_r,l)$
in the new relations. Beside this, the following natural modification
of formula (3.11) holds.

If the decomposition of a diagram
$\gamma=\gamma_A\in\Gamma(k_j|j\in A)$ to connected components
consists of the diagrams
$\gamma_{A_1}$,\dots, $\gamma_{A_u}$, where $A_1,\dots,A_u$ is
the partition of the set $A$ we have to apply to get the desired
decomposition to connected components, then the above defined
function $F_\gamma(x_n,\,n\in N_1(\gamma))$ satisfies the identity
$$
F_{\gamma_A}(x_n,\,n\in N_1(\gamma_A))=\prod_{t=1}^u
F_{\gamma_{A_t}}(x_n,\,n\in N_1(\gamma_{A_t})),
\tag4.3
$$
where $N_1(\gamma_{A_t})$ denotes the set of labels of the
open vertices in $\gamma_{A_t}$, and $F_{\gamma_{A_t}}$ is the
function defined in (4.2) with the diagram $\gamma=\gamma_{A_t}$.
Beside this, the sets of labels $N_1(\gamma_{A_t})$ are disjoint
for different indices~$t$. Here I make the convention that the
label of a vertex in the restriction $\gamma_B$ of a diagram
$\gamma_A$, $B\subset A$, agrees with its original label in the 
diagram $\gamma_A$. Formula (4.3) can be checked similarly
to~(3.8) with the help of (the modified version of) the recursive
relations~(3.8) and~(3.9).

Given a closed diagram $\gamma$, its restriction $\gamma^r$ to its
first $r$ rows may be not a closed diagram. But it is also a
diagram, and the above mentioned results can be applied for it. In
particular, it has a decomposition to connected components,
and the function $F_{\gamma^r}$ defined by means of some functions
$f_j$, $1\le j\le r$ can be factorized. The next result, called
{\it A More Detailed Version of the Basic Estimate}\/ contains an
estimate for the terms in this factorization. The Basic Estimate is
a part of this result for the parameter $r=m$.

\medskip\noindent
{\bf  A More Detailed Version of the Basic Estimate.}
{\it Let us consider a\/ {\rm connected}, closed diagram
$\gamma\in\bar\Gamma(k_1,\dots,k_m)$, $m\ge2$, and some functions
$f_j$ of $k_j$ variables on a measure space $(X,\Cal X,\mu)$,
$1\le j\le m$, which satisfy the inequality $V_s(f_j)\le R^{s-1}$
with some $0\le R\le1$ for all $1\le j\le m$ and $1\le s\le k_j$.
Let $\gamma^r$ denote the restriction of the diagram $\gamma$ to
its first $r$ rows and let $F_{\gamma^r}$ be the function defined
in formulas~(4.1) and~(4.2) with the help of the diagram $\gamma=\gamma^r$
and the functions $f_j$, $1\le j\le r$. Take the decomposition of
the diagram $\gamma^r$ to the union of disjoint connected
diagrams $\gamma_{A^r_t}$ defined with the help of an appropriate
partition $A^r_1,\dots,A^r_{u(r)}$ of the set $\{1,\dots,r\}$
with some number $u(r)$. Consider the factorization (4.3) of the
function $F_{\gamma^r}$ to
$$
F_{\gamma^r}(x_n,\; n\in N_1(\gamma^r))=\prod_{t=1}^{u(r)}
F_{\gamma_{A^r_t}}(x_n,\,n\in N_1(\gamma_{A^r_t})).
$$
For $1\le r\le m-1$ all diagrams $\gamma_{A^r_t}$, $1\le t\le u(r)$,
contain at least one open vertex, i.e. $N_1(\gamma_{A^r_t})\ge1$,
and
$$
|V_s(F_{\gamma_{A^r_t}}(x_n,\,n\in N_1(\gamma_{A^r_t})))|\le
R^{|A^r_t|+s-2} \quad\text{for all } 1\le t\le u(r)
\text{ and } 1\le s\le |N_1(\gamma_{A^r_t})|,
\tag4.4
$$
where $|A^r_t|$ is the cardinality of the set $A^r_t$, i.e. it equals
the number of rows in the diagram~$A^r_t$, and $N_1(\gamma_{A^r_t})$
denotes the set of indices of the free vertices in $\gamma_{A^r_t}$.

For $r=m$ the decomposition of $\gamma^m=\gamma$ to connected
components consists of the (closed) diagram $\gamma$ itself, and
$$
|F_\gamma|=|F_{\gamma^m}|\le R^{m-2}. \tag4.5
$$
}

\medskip
The proof of the {\it More Detailed Version of the Basic Estimate}\/
is based on the following result called the Main Inequality.

\medskip\noindent
{\bf The Main Inequality.} {\it Let
$f(x_1,\dots,x_m,v_{m+n+1},\dots,v_{m+n+q})$
and
$$
g(x_{m+1},\dots,x_{m+n},v_{m+n+1},\dots,v_{m+n+q})
$$
be two square integrable functions with $m+q$ and $n+q$ variables on a 
measure space $(X,\Cal X,\mu)$, and define the function
$$
\aligned
F(x_1,\dots,x_{m+n})
&=\int f(x_1,\dots,x_m,v_{m+n+1},\dots,v_{m+n+q})\\
&\qquad g(x_{m+1},\dots,x_{m+n},v_{m+n+1},\dots,v_{m+n+q})
\prod_{u=m+n+1}^{m+n+q}\mu(\,dv_u).
\endaligned \tag4.6
$$
Let $m+n\ge1$, $q\ge1$, and let the functions $f$ and
$g$ satisfy the relation $V_s(f)\le D_1R^{s-2}$ and
$V_s(g)\le D_2R^{s-2}$ with some $D_1>0$, $D_2>0$ and $0\le R\le 1$ 
for all $1\le s\le m+q$ and $1\le s\le n+q$ respectively. Then the
function $F$ satisfies the inequality
$$
V_s(F(x_1,\dots,x_{m+n}))\le D_1D_2 R^{s-2}
\quad\text{for all \ } 1\le s\le m+n. \tag4.7
$$
}

First I prove the More Detailed Version of the Basic Estimate with
the help of the Main Inequality.

\medskip\noindent
{\it Proof of the More Detailed Version of the Basic Estimate with
the help of the Main Inequality.} Let us first observe that all
components $\gamma_{A^r_t}$, $1\le t\le u(r)$, of the diagram
$\gamma^r$ contain at least one open vertex for $r\le m-1$.
Otherwise the diagram $\gamma$ would not satisfy the condition
of connectedness with $B=A^r_t$. Relation (4.4) clearly holds for
$r=1$. It will be proved by induction with respect to $r$ that it
holds for all $1\le r\le m-1$.

Let us assume that relation (4.4) holds for $r-1$ and let us prove
it for $r$ if $r\le m-1$. Some of the connected components 
$\gamma_{A^{r-1}_t}$ of $\gamma^{r-1}$ may have a vertex connected 
with a vertex of the $r$-th row $\{(r,1),\dots,(r,k_r)\}$ of the 
diagram~$\gamma$. Let us make such an enumeration of the connected 
components of $\gamma^{r-1}$ in which there is a constant 
$\bar u(r)$ such that the components $\gamma_{A^{r-1}_t}$ with 
$1\le t\le\bar u(r)$ have a vertex connected with a vertex of the 
$r$-th row of~$\gamma$, and the components $\gamma_{A^{r-1}_t}$, 
$\bar u(r)<t\le u(r-1)$, have no such vertices. To simplify the 
following discussion let us introduce the notation 
$\gamma_{A^{r-1}_0}=\{(r,1),\dots,(r,k_r)\}$ for the $r$-th row of 
the diagram~$\gamma$. Then the connected components of $\gamma^r$ 
are the diagrams 
$\gamma_{B^r}=\bigcupp_{t=0}^{\bar u(r)}\gamma_{A^{r-1}_t}$
and $\gamma_{A^{r-1}_t}$, $\bar u(r)<t\le u(r-1)$. The latter
diagrams satisfy relation (4.4) by induction, so it is enough to
show that
$$
|V_s(F_{\gamma_{B^r}})|\le R^{|B^r|+s-2} \quad\text{for all }
1\le s\le |N_1(\gamma_{B^r})|, \tag4.8
$$
where $|N_1(\gamma_{B^r})|$ denotes the number of open vertices in
the diagram $\gamma_{B^r}$. It is possible that no vertex of
the diagram $\gamma^{r-1}$ is connected with a vertex from the
$r$-th row of $\gamma$. In this case $\bar u(r)=0$, $B^r=A^{r-1}_0$,
and relation~(4.8) clearly holds.

To prove relation (4.8) let us introduce the following notations.
Put $B^r_j=\bigcupp_{t=0}^j A^{r-1}_t$ for all $0\le j\le \bar u(r)$.
With such a notation $B^r_{\bar u(r)}=B^r$ for $j=\bar u(r)$. I shall 
show that
$$
|V_s(F_{\gamma_{B^r_j}}(x_n,\; n\in N_1(\gamma_{B^r_j})|
\le R^{|B^r_j|+s-2} \quad\text{for all } 1\le j\le \bar u(r)
\text { and } 1\le s\le |N_1(\gamma_{B^r_j})|. \tag4.9
$$
Relation (4.9) for $j=\bar u(r)$ implies relation (4.8). Relation~(4.9)
holds for $j=0$, and it will be proved for a general parameter $j$,
$1\le j\le \bar u(r)$, by induction.

For this goal I write the following recursive relation for the
functions $F_{\gamma_{B^r_j}}$
$$
\aligned
F_{\gamma_{B^r_j}}(x_n,\; n\in N_1(\gamma_{B^r_j}))&=\int
F_{\gamma_{A^{r-1}_j}}(x_n,\; n\in N_1(\gamma_{A^{r-1}_j}))\\
&\qquad F_{\gamma_{B^r_{j-1}}}(x_n,\; n\in N_1(\gamma_{B^r_{j-1}}))
\prod_{n\in N_1(\gamma_{A^{r-1}_j})
\cap N_1(\gamma_{B^r_{j-1}})}\mu(\,dx_n)
\endaligned \tag4.10
$$
for all $1\le j\le \bar u(r)$. I show that relation (4.9) follows
from relation~(4.10) and the Main Inequality. Indeed, let us
apply the Main Inequality with the functions
$f=F_{\gamma_{A^{r-1}_j}}(x_n,\; n\in N_1(\gamma_{A^{r-1}_j}))$
and $g=F_{\gamma_{B^r_{j-1}}}(x_n,\;n\in N_1(\gamma_{B^r_{j-1}}))$.
(More precisely, we apply an equivalent version of the Main
Inequality where the indices of the variables of the functions $f$
and $g$ may be different, and the variables by which we integrate
and by which we do not integrate may be listed in an arbitrary
order.) By our inductive hypothesis these functions $f$ and $g$
satisfy the inequalities
$V_s(f)\le D_1R^{s-2}$ and $V_s(g)\le D_2R^{s-2}$ with
$D_1=R^{|A^{r-1}_j|}$ and $D_2=R^{|B^r_{j-1}|}$. Beside this,
$|N_1(\gamma_{A^{r-1}_j})\cap N_1(\gamma_{B^r_{j-1}})|=
|N_1(\gamma_{A^{r-1}_j})\cap N_1(\gamma_{A^{r-1}_0})|\ge1$, because
there is an edge connecting a vertex of $\gamma_{A^{r-1}_j}$ with a
vertex of the $r$-th row $\gamma_{A^{r-1}_0}$ of the diagram
$\gamma$. This inequality corresponds to the condition $q\ge1$
in the Main Inequality, where the number $q$ is the multiplicity of 
the integral in formula~(4.6). Beside this, the diagram
$\gamma_{A^{r-1}_j}\cup \gamma_{B^r_{j-1}}$ has an open vertex
because of the connectedness of the diagram $\gamma$. This
corresponds to the condition $m+n\ge1$ in the Main Inequality.

The above considerations show that the Main Inequality can be
applied in the present case. It yields that
$|V_s(F_{\gamma_{B^r_j}}|\le R^{|A^{r-1}_j|+|B^r_{j-1}|}R^{s-2}
=R^{|B^r_j|+s-2}$, and this is what we had to prove.

The proof of relation (4.5) is similar. In the proof the
above decomposition of $\gamma^{r-1}$ is applied to the connected
components for $r=m$. For the parameter $r=m$ all components
$\gamma_{A^{m-1}_t}$, $1\le t\le u(m-1)$ have a vertex which is
connected with a vertex of the $m$-th row of the diagram $\gamma$.
(The $m$-th row of $\gamma$ will be sometimes denoted by
$\gamma_0^{m-1}$.) Thus $\bar u(m)=u(m-1)$, the connected components
$\gamma_{A^{m-1}_j}$, $1\le j\le u(m-1)$, can be listed in an
arbitrary order, and the diagrams $\gamma_{B^m_j}$, $1\le j\le u(m-1)$
can be defined similarly to the case $r<m$. Beside this,
relation~(4.9) can be proved for $r=m$ and $j\le u(m-1)-1$,
similarly as it was proved for $r<m$. The main difference between
the cases $r<m$ and $r=m$ is that in the latter case the proof of
relation~(4.9) works only for $j\le u(m-1)-1$, but not for
$j=u(m-1)$. In the case $j=u(m-1)$ the Main Inequality
cannot be applied in the proof, because $\gamma_{B^m_{u(m)}}=\gamma$
is a diagram without open vertices, and this property does not allow
the application of the Main Inequality in the case $r=m$, $j=u(m-1)$.

On the other hand the identity
$$ \allowdisplaybreaks
\align
F_\gamma=&\int
F_{\gamma_{B^m_{u(m-1)-1}}} (x_n,\;n\in
N_1(\gamma_{A^{m-1}_{u(m-1)} })) \\
&\qquad F_{\gamma_{A^{m-1}_{u(m-1)}}}(x_n,\;n\in
N_1(\gamma_{A^{m-1}_{u(m-1)}}))
\prod_{n\in N_1(\gamma_{A^{m-1}_{u(m-1)}})}\mu(\,dx_n) \tag4.11
\endalign
$$
holds, and the function integrated in (4.11) is the product of two
terms whose $L_2$-norm can be well bounded. Namely, since the 
$L_2$-norm of a function $f$ of several variables (defined on the 
measure space $(X,\Cal X,\mu)$) equals $V_1(f)$, relation (4.9) with 
the choice $r=m$ and $j=u(m-1)-1$ together with formula (4.4) for 
$r=m-1$ yield, with the parameter $s=1$, the bound 
$R^{|B^m_{u(m-1)-1}|-1}$
and $R^{|A^{m-1}_{u(m-1)}|-1}$ for the $L_2$-norm of these terms.
Hence relation (4.11) and the Schwarz inequality imply that
$$
|F_\gamma|\le R^{|B^m_{u(m-1)-1}|+A^{m-1}_{u(m-1)}|-2}=R^{m-2},
$$
as it was stated.

\medskip
It remained to prove the Main Inequality.

\medskip\noindent
{\it Proof of the Main Inequality.}\/ To formulate the inequality
we have to prove first some notation will be introduced. A 
partition of the set $\{1,\dots,m+n\}$ will be introduced which
tells the indices of the functions $u_j(\cdot)$ we shall work with.
This partition will consist of $s=s_1+s_2+s_3$ elements, where $s_1$
is the number of those sets in this partition which have an element
in both sets $\{1,\dots,m\}$ and $\{m+1,\dots,m+n\}$, $s_2$ and $s_3$ 
are the number of the sets in the partition which are contained in the
set $\{1,\dots,m\}$ and $\{m+1,\dots,m+n\}$ respectively. In the
first step of the proof it will be shown that we can restrict our
attention to the case $s_1=s$, $s_2=s_3=0$.

The following notation will be used. Let us fix a partition
$A_j=A_j^{(1)}\cup A_j^{(2)}$, $1\le j\le s_1$, $B_k$, $1\le k\le s_2$,
and $C_l$, $1\le l\le s_3$, $s_1+s_2+s_3=s$ of the set
$\{1,\dots,m+n\}$ such that all sets $A^{(1)}_j$, $A^{(2)}_j$, $B_k$
and $C_l$ are non-empty, and $A^{(1)}_j, B_k\subset\{1,\dots,m\}$,
$A^{(2)}_j, C_l\subset \{m+1,\dots,m+n\}$. Define some functions
$u_j(x_p,\;p\in A_j)$, $1\le j\le s_1$, $v_k(x_r,\; r\in B_k)$,
$1\le k\le s_2$ and $w_l(x_t,\;t\in C_l)$, $1\le l\le s_3$ such
that $\int u^2_j(x_p,\;p\in A_j)\prodd_{p\in A_j}\mu(\,dx_p)\le1$,
$\int v^2_k(x_r,\; r\in B_k)\prodd_{r\in B_k}\mu(\,dx_r)\le1$,
and $\int w^2_l(x_t,\;t\in C_l)\prodd_{t\in C_l}\mu(\,dx_t)\le1$.
It has to be shown that for all such partitions and functions
$u_j(\cdot)$,  $v_k(\cdot)$ and $w_l(\cdot)$ the functions $f$ and 
$g$ satisfy the inequality
$$
\aligned
&\int f(x_1,\dots,x_m,v_{m+n+1},\dots,v_{m+n+q})
g(x_{m+1},\dots,x_{m+n},v_{m+n+1},\dots,v_{m+n+q}) \\
&\qquad \prod_{j=1}^{s_1} u_j(x_p,\;p\in A_j)
\prod_{k=1}^{s_2} v_k(x_r,\; r\in B_k)
\prod_{l=1}^{s_3} w_l(x_t,\;t\in C_l)\\
&\qquad\qquad \prod_{i=1}^{m+n}\mu(\,dx_i)
\prod_{u=m+n+1}^{m+n+q}\mu(\,dv_u) \le D_1D_2 R^{s-2}.
\endaligned \tag4.12
$$
First we reduce this inequality to the case $s_2=s_3=0$, i.e. to the
case when the partition of the set $\{1,\dots,n+m\}$ consists only of 
such sets $A_j$ which have a non-empty intersection with both sets 
$\{1,\dots,n\}$ and $\{n+1,\dots,n+m\}$. For this goal we define the
functions
$$
\aligned
&\bar f(x_j,\;j\in \{1,\dots,m\}\setminus B,
v_{m+n+1},\dots,v_{m+n+q})\\
&\qquad =\int f(x_1,\dots,x_m,v_{m+n+1},\dots,v_{m+n+q})
\prod_{k=1}^{s_2} v_k(x_r,\; r\in B_k)
\prod_{r\in B}\mu(\,dx_r),
\endaligned \tag4.13
$$
$$
\aligned
&\bar g(x_j,\;j\in \{m+1,\dots,m+n\}\setminus C,
v_{m+n+1},\dots,v_{m+n+q})\\
&\qquad =\int g(x_{m+1},\dots,x_{m+n},v_{m+n+1},\dots,v_{m+n+q})
\prod_{l=1}^{s_3} w_l(x_t,\;t\in C_l)
\prod_{t\in C}\mu(\,dx_t)
\endaligned \tag4.14
$$
with $B=\bigcupp_{k=1}^{s_2} B_k$ and $C=\bigcupp_{l=1}^{s_2} C_l$
together with 
$$
\aligned
&\bar F(x_j,\; j\in\{1,\dots,m+n\}\setminus (B\cup C) \\
&\qquad=\int
\bar f(x_j,\;j\in \{1,\dots,m\}\setminus B,
v_{m+n+1},\dots,v_{m+n+q}) \\
&\qquad\qquad\quad
\bar g(x_j,\;j\in \{m+1,\dots,m+n\}\setminus C,
v_{m+n+1},\dots,v_{m+n+q}) \prod_{u=m+n+1}^{m+n+q}\mu(\,dv_u).
\endaligned \tag4.15
$$
With this notation inequality (4.12) can be rewritten
as
$$ \allowdisplaybreaks
\align
\int &\bar F(x_i,\; i\in\{1,\dots,m+n\}\setminus (B\cup C)
\prod_{j=1}^{s_1} u_j(x_q,\;q\in A_j) \\
&\qquad \prod_{i\in\{1,\dots,m+n\}\setminus (B\cup C)} \mu(\,dx_i)
\le D_1D_2 R^{s_1+s_2+s_3-2}. \tag4.16
\endalign 
$$
Beside this, it is not difficult to check that the functions $\bar f$
and $\bar g$ satisfy the inequalities $V_s(\bar f)\le \bar D_1R^{s-2}$
for all $1\le s\le n+q-|B|$ and $V_s(\bar g)\le \bar D_2R^{s-2}$
for all $1\le s\le m+q-|C|$ with $\bar D_1=D_1R^{s_2}$ and
$\bar D_2=D_2R^{s_3}$ respectively. For the parameter $s=1$ these
inequalities yield the bounds $\bar D_1R^{-1}$ and $\bar D_2R^{-1}$
for the $L_2$-norm of the functions $\bar f$ and $\bar g$ 
respectively. They imply together with relation (4.15) and
the Schwarz inequality that relation (4.16) holds in the special case
$s_1=0$, i.e. when there is no set of the type $A_j$ in the partition 
we consider. This special case had to be considered separately, 
because in this case $\bar F$ is a function of zero variables, i.e. 
it is a constant. This means that in the reduced model we want to
consider the condition $m+n\ge1$ is violated in this special case.

By the above observation it is enough to prove relation (4.16) in 
the case $s_1\ge1$. This enables us to reduce the proof of 
relation~(4.12) to the case $s_2=s_3=0$, i.e. to the case when all 
elements of the partition is such a set as the sets $A_j$. We get 
this reduction by working with the functions $\bar f$, $\bar g$ and
$\bar F$ defined in formulas (4.13), (4.14) and (4.15) instead of
the original functions $f$, $g$ and $F$ in the proof of the Main
Inequality, and by observing that these functions also satisfy its
conditions (with an appropriate reindexation of the variables in
these functions). Hence the Main inequality follows from the
following reduced form of relation (4.12):
$$
\aligned
I&=\int f(x_1,\dots,x_m,v_{m+n+1},\dots,v_{m+n+q})
g(x_{m+1},\dots,x_{m+n},v_{n+m+1},\dots,v_{m+n+q}) \\
&\qquad \prod_{j=1}^s u_j(x_p,\;p\in A_j)
\prod_{i=1}^{m+n}\mu(\,dx_i)
\prod_{u=m+n+1}^{m+n+q}\mu(\,dv_u) \le D_1D_2 R^{s-2}
\endaligned \tag4.17
$$
for a partition $A_j=A^{(1)}_j\cup A_2^{(j)}$, $1\le j\le s$, of 
the set $\{1,\dots,m+n\}$. To prove inequality~(4.17) let us 
introduce the functions
$$
U_j(x_p,\;p\in A_j^{(2)})=\[\int u_j^2(x_p,\;p\in A_j)
\prodd_{r\in A_j^{(1)}}\mu(\,dx_r)\]^{1/2}
$$
and
$$
\bar u_j(x_p,\;p\in A_j^{(1)}|x_p,\; p\in A_j^{(2)})=
\frac{u_j(x_p,\;p\in A_j)}{U_j(x_p,\;p\in A_j^{(2)})}
$$
for all $1\le j\le s$ together with the functions
$$
\align
&G(x_{m+1},\dots,x_{m+n})\\
&\qquad =\[\int
g^2(x_{m+1},\dots,x_{m+n},v_{m+n+1},\dots,v_{m+n+q})
\prod_{u=m+n+1}^{m+n+q}\mu(\,dv_u)\]^{1/2}
\endalign
$$
and
$$
\bar g(v_{m+n+1},\dots,v_{m+n+q}|x_{m+1},\dots,x_{m+n})=
\frac{g(x_{m+1},\dots,x_{m+n},v_{m+n+1},\dots,v_{m+n+q})}
{G(x_{m+1},\dots,x_{m+n})}.
$$
Observe that
$$
\int U_j^2(x_p,\;p\in A_j^{(2)})
\prodd_{p\in A_j^{(2)}}\mu(\,dx_p)\le 1, \tag4.18
$$
$$
\int \bar u^2_j(x_p,\;p\in A_j^{(1)}|x_p,\; p\in A_j^{(2)})
\prod_{p\in A_j^{(1)}}\mu(\,dx_p)=1 \quad\text{ for all }
x=\{x_p,\; p\in A_j^{(2)}\} \tag4.19
$$
for all $1\le j\le s$, because the $L_2$-norm of the functions
$u_j(\cdot)$ are less than~1. Similarly,
$$
\int G^2(x_{m+1},\dots,x_{m+n})\prod_{j=m+1}^{m+n}\mu(\,dx_j)
\le D_2^2R^{-2}, \tag4.20
$$
since the $L_2$-norm of the function $G$ equals the $L_2$-norm of
$g$ which is $V_1(g)$, and it is bounded by $D_2R^{-1}$. Beside this,
$$
\aligned
\int &\bar g^2(v_{m+n+1},\dots,v_{m+n+q}|x_{m+1},\dots,x_{m+n})
\prod_{u=m+n+1}^{m+n+q}\mu(\,dv_u)=1 \\
&\qquad \text{for all } x=(x_{m+1},\dots,x_{m+n})
\endaligned \tag4.21
$$

The expression $I$ in formula (4.17) can be rewritten with the help 
of the identities $u_j(\cdot)=\bar u_j(\cdot|\cdot)U_j(\cdot)$, 
$1\le j\le s$ and $g(\cdot)=\bar g(\cdot|\cdot)G(\cdot)$ as
$$
I=\int Z(x_r,\;m+1\le r\le  m+n)
\prod _{j=1}^s U_j(x_p,\;p\in A_j^{(2)})
G(x_{m+1},\dots,x_{m+n})\prod_{i=m+1}^{m+n}\mu(\,dx_i). \tag4.22
$$
with the help of the function
$$
\align
Z(x_r,\;m+1\le r\le  m+n)
&=\int f(x_1,\dots,x_m,v_{m+n+1},\dots,v_{m+n+q})\\
&\qquad \prod_{j=1}^s
\bar u_j(x_p,\;p\in A_j^{(1)}|x_p,\; p\in A_j^{(2)})\\
&\!\!\!\!\!\!\!\!\!\!\!\!\!\!\!\!\!\!\!\!\!\!\!\!
\bar g(v_{m+n+1},\dots,v_{m+n+q}|x_{m+1},\dots,x_{m+n})
\prod_{i=1}^m\mu(\,dx_i)\prod_{u=m+n+1}^{m+n+q}\mu(\,dv_u).
\endalign
$$

The function $Z$ satisfies the inequality
$$
|Z(x_r,\;m+1\le r\le m+n)|\le D_1R^{s-1} \quad \text{for all }
x=(x_r,\;m+1\le r\le m+n), 
$$
because of the inequality $V_{s+1}(f)\le D_1 R^{s-1}$ (we consider
the partition of the set $\{1,\dots,m,m+n+1,\dots,m+n+q\}$ 
consisting of the sets $A^{(1)}_j$, $1\le j\le s$, and the set 
$\{m+n+1,\dots,m+n+q\}$) 
and the bounds (4.19) and (4.21) about the $L_2$-norm of the 
functions $\bar u_j(\cdot|\cdot)$ and $\bar g(\cdot|\cdot)$.

Beside this, the $L_2$-norms of the functions
$$
\prodd _{j=1}^s U_j(x_p,\;p\in A_j^{(2)}) \quad \text{and}
\quad G(x_{m+1},\dots,x_{m+n})
$$
are bounded by 1 and $D_2R^{-1}$ respectively by relations~(4.18) 
and~(4.20).  These inequalities together with relation~(4.22) and 
the Schwarz inequality imply that $I\le D_1D_2R^{s-2}$, i.e. 
relation~(4.17) holds. The proof of the Main Inequality is 
completed.

\beginsection 5. On Lata{\l}a's conjecture.

In this section I discuss Lata{\l}a's conjecture. I show,
by working out the details of the arguments leading to this 
conjecture that it is equivalent to an estimate about the 
expected value of the supremum of certain random multilinear 
forms. The original form of this conjecture contains an 
estimate about the moments of certain Gaussian polynomials. The 
would-be proof applies an inductive argument with respect to 
the order of the polynomials we consider. It is known that the 
conjecture holds for polynomials of order~$2$. I shall consider 
polynomials of order~3 whose study also reveals very much about
the general situation. First I formulate Lata{\l}a's conjecture
for polynomials of order three in an explicit form.

\medskip\noindent
{\bf Lata{\l}a's conjecture for Gaussian random polynomials of
order~3.} {\it Let us fix a large positive integer $M$ and
consider a random polynomial of the form
$$
Z=\sum_{ i,j,k} a(i,j,k)\xi_i\eta_j\zeta_k, \tag5.1
$$
where all random variables $\xi_i$, $\eta_j$ and $\zeta_k$ have
standard normal distribution, and they are independent of each
other. Let the coefficients $a(i,j,k)$ of the polynomial $Z$ in
formula (5.1) satisfy the following inequalities depending on
the fixed parameter~$M$:
$$
\sum_{i,j,k} a(i,j,k)u(i,j,k)\le 1 \quad\text{if }
\sum_{i,j,k}u^2(i,j,k)\le1, \tag5.2
$$
$$
\aligned
\sum_{i,j,k} a(i,j,k)u(i,j)v(k)&\le M^{-1/2} \quad\text{if }
\sum_{i,j}u^2(i,j)\le1 \quad \text{and } \sum_{k} v^2(k)\le1,\\
\sum_{i,j,k} a(i,j,k)u(i,k)v(j)&\le M^{-1/2} \quad\text{if }
\sum_{i,k}u^2(i,k)\le1 \quad \text{and } \sum_{j}v^2(j)\le1,\\
\sum_{i,j,k} a(i,j,k)u(j,k)v(i)&\le M^{-1/2} \quad\text{if }
\sum_{j,k}u^2(j,k)\le1 \quad \text{and } \sum_{i}v^2(i)\le1,
\endaligned \tag5.3
$$
and
$$
\align
\sum_{i,j,k} &a(i,j,k)u(i)v(j)w(k)\le M^{-1} \\
&\qquad\text{if} \quad \sum_{i}u^2(i)\le1, \quad \sum_{j} v^2(j)\le1
\quad \text{and } \sum_{k} w^2(k)\le1. \tag5.4
\endalign
$$
Then the random polynomial $Z$ satisfies the inequality
$$
EZ^{2M}\le C^MM^M \tag5.5
$$
with some universal constant $C>0$.}

\medskip
In the calculation of  $EZ^{2M}$ it is useful to consider first its
conditional expectation under the condition that the value of all
random variables $\xi_i$ are prescribed. This conditional
expectation has a simple form which can be well bounded because of
the independence of the variables $\xi_i$, $\eta_j$ and $\zeta_k$.
We get
$$
E(Z^{2M}|\xi_i=x_i)=E\(\sum_{ i,j,k} a(i,j,k)x_i\eta_j\zeta_k\)^{2M}
=E\(\summ_{j,k}A(j,k|x)\eta_j\zeta_k\)^{2M}, \tag5.6
$$
where
$$
A_i(j,k|x)=A_i(j,k|x_1,x_2,\dots)=\summ_i a(i,j,k)x_i.
$$
The moment estimates known for Gaussian polynomials of order two
enable us to bound the expression in formula~(5.6). These estimates
depend on the Hilbert--Schmidt norm $D_1(x)$ and usual norm $D_2(x)$ 
of the matrix $A(j,k|x)$ appearing in formula~(5.6). To get a formula
more appropriate for our investigations let us give the value of these
quantities by means of the following variational principle.
$$
D_1(x)=\supp_{v(j,k)\colon\; \summ v^2(j,k)\le1}
\summ_{j,k}A(j,k|x)v(j,k),
$$
and
$$
D_2(x)=\supp_{(v(j),\,w(k))\colon\; \summ v^2(j)\le1,\;\summ w^2(k)\le 1}
\summ_{j,k}A(j,k|x)v(j)w(k).
$$
In such a way we get the following estimate.
$$
\aligned
&E(Z^{2M}|\xi_i=x_i)=E\(\summ_{j,k}A(j,k|x)\eta_j\zeta_k\)^{2M}\\
&\le C^M M^M\(\sup\Sb v(j,k)\\ \summ v^2(j,k)\le1\endSb
\summ_{j,k}A(j,k|x)v(j,k)\)^{2M} \\
&\qquad + C^MM^{2M}\(\sup\Sb v(j),\;w(k)\\ \summ v^2(j)\le1,\;\;
\summ w^2(k)\le 1\endSb
\summ_{j,k}A(j,k|x)v(j)w(k)\)^{2M}.
\endaligned \tag5.7
$$

Taking expectation in inequality (5.7) we get that
$$
\aligned
E(Z^{2M})
&\le C^M M^M E\(\sup\Sb u(j,k)\\ \summ u^2(j,k)\le1\endSb
\summ_{j,k}A(j,k|\xi_i)v(j,k)\)^{2M}\\
&\qquad+ C^M M^{2M} E\(\sup\Sb v(j),\;w(k)\\ \summ v^2(j)\le1,\;\;
\summ w^2(k)\le 1\endSb
\summ_{j,k}A(j,k|\xi_i)v(j)w(k)\)^{2M}.
\endaligned \tag5.8
$$
The last inequality can be rewritten in the form
$$
\aligned
E(Z^{2M})
&\le C^M M^M E\(\sup\Sb u(j,k)\\ \summ u^2(j,k)\le1\endSb
\summ_{i}B_{i,1}(u(j,k))\xi_i\)^{2M}\\
&\qquad+ C^M M^{2M} E\(\sup\Sb v(j),\;w(k)\\ \summ v^2(j)\le1,\;\;
\summ w^2(k)\le 1\endSb
\summ_{i}B_{i,2}(v(j),\,w(k))\xi_i\)^{2M}.
\endaligned 
$$
with
$$
B_{i,1}(u(j,k),\; j,k=1,2,\dots,)=\summ_{j,k} a(i,j,k)u(j,k),
$$
and
$$
B_{i,2}(v(j),\, w(k),\;j,k=1,2,\dots)=\summ_{j,k} a(i,j,k)v(j)w(k),
$$
or by introducing the notations $u=(u(j,k),\;j,k=1,2,\dots)$,
$v=(v(j),\;j=1,2,\dots)$ and $w=(v(k),\;k=1,2,\dots)$ together
with the Gaussian random variables
$$
\aligned
X(u)&=\sum_{i,j,k} a(i,j,k) u(j,k)\xi_i
=\sum_{i} B_{i,1}(u(j,k),\;j,k=1,2,\dots)\xi_i, \\
Y(v,w)&=\sum_{i,j,k} a(i,j,k) v(j)w(k)\xi_i
=\sum_{i} B_{i,2}(v(j),\,w(k),\;j,k=1,2,\dots)\xi_i, \
\endaligned \tag5.9
$$
this inequality can be written in the form
$$
\aligned
E(Z^{2M})
&\le C^M M^M E\(\sup\Sb u(j,k)\\ \summ u^2(j,k)\le1\endSb
X(u)\)^{2M} \\
&\qquad + C^M M^{2M} E\(\sup\Sb v(j),\;w(k)\\ \summ v^2(j)\le1,\;\;
\summ w^2(k)\le 1\endSb Y(v,w)\)^{2M}.
\endaligned \tag5.10
$$

The right-hand side of (5.10) can be bounded by means of some
concentration theorem type inequalities about the supremum of a
Gaussian process. Ledoux has a result about the supremum of
Gaussian  processes (Theorem 7.1 in the book~[7] {\it The
Concentration of Measure Phenomenon}) which states that the
supremum of a Gaussian process $U(t)$, $EU(t)=0$, $t\in T$,
takes a value larger than $E\supp_{t\in T}U(t)$ with relatively
small probability. More explicitly, it states that
$P\(\supp_{t\in T} U(t)\ge E\supp_{t\in T} U(t)+x\)\le
C_1e^{-C_2x^2/\lambda}$, where $\lambda=\supp_{t\in T}EU^2(t)$.
Some calculation with the help of this inequality yields the
estimates
$$
\align
&E\(\sup\Sb u(j,k)\\ \summ u^2(j,k)\le1\endSb X(u)\)^{2M}
\le D^M
\(E\left|\sup\Sb u(j,k)\\ \summ u^2(j,k)\le1\endSb X(u)\right|\)^{2M}\\
&\qquad\qquad
+D^M\supp\Sb u(j,k) \\ \sum u^2(j,k)\le 1\endSb  EX(u)^{2M}\\
&\le D^M \(E\left|\sup\Sb u(j,k)\\
\summ u^2(j,k)\le1\endSb X(u)\right|\)^{2M}
+D^M M^M\supp\Sb u(j,k) \\ \sum u^2(j,k)\le 1\endSb \(EX(u)^2\)^M
\endalign
$$
and
$$
\align
&E\(\sup\Sb v(j),\;w(k)\\ \summ v^2(j)\le1,\;\;
\summ w^2(k)\le 1\endSb Y(v,w)\)^{2M} \\
&\qquad \le D^M
\(E\left|\sup\Sb v(j),\;w(k)\\
\summ v^2(j)\le1,\;\;\sum w^2(k)\le1\endSb Y(v,w)\right|\)^{2M} \\
&\qquad\qquad
+D^M \supp\Sb v(j),\;w(k) \\ \sum v^2(j)\le1,\;\;w^2(k)\le 1\endSb
EY(v,w)^{2M}\\
&\le D^M\( E\left|\sup\Sb v(j),\;w(k)\\ \summ v^2(j)\le1,\;\;
\sum w(k)^2\le1\endSb Y(v,w)\right|\)^{2M} \\
&\qquad\qquad +D^M M^M
\supp\Sb v(j),\;w(k) \\ \sum v^2(j)\le1,\;\;\sum w^2(k)\le 1\endSb
\(EY(v,w)^2\)^M.
\endalign
$$

The content of the above inequalities is that to get a good estimate
on the high moments of the supremum of a Gaussian process with
expectation zero it is enough to have a good estimate on the
expectation of the absolute value of this supremum and on the moments
of the single random variables in this stochastic process. The
latter terms can be  expressed by means of the variance of these
random variables.

A relatively simple calculation by means of the Schwarz inequality
shows that under conditions (5.3) and (5.4) the inequalities
$EX(u)^2\le M^{-1}$, and $EY(v,w)^2\le M^{-2}$ hold if
$\sum u(j,k)^2\le1$, $\sum v^2(j)\le1$, and $\sum w^2(k)\le1$.
Some calculation also shows that under the condition (5.2) the
inequality
$$
E\left|\sup\Sb u(j,k)\\ \summ u^2(j,k)\le1\endSb X(u)\right|\le
\(E\left(\sup\Sb u(j,k)\\ \summ u^2(j,k)\le1\endSb X(u)\right)^2\)^{1/2}
\le C \tag5.11
$$
holds. The second term in~(5.11) can be well estimated, because the
supremum of the random variables $X(u)$ can be explicitly
calculated for all fixed random vectors $\xi_i$, $i=1,2,\dots$, and
after this, the middle term in (5.11) can be well bounded with the
help of relation (5.2) because of the orthogonality of the random
variables $\xi_i$.

Because of the above inequalities to show that relation (5.5) holds 
under conditions (5.2), (5.3) and~(5.4) it would be sufficient to 
prove the inequality
$$
E\left|\sup\Sb v(j),\;w(k)\\ \summ v^2(j)\le1,\;\;
\sum w(k)^2\le1\endSb Y(v,w)\right|\le \frac C{\sqrt M}. \tag5.12
$$
under the above conditions with the random variables $Y(v,w)$
introduced in (5.9) and some universal constant $C<\infty$. This
would mean that Lata{\l}a's conjecture holds for Gaussian
polynomials of order~3.

A more careful analysis would even show that the validity of
relation~(5.12) under conditions (5.2), (5.3) and~(5.4) is
equivalent to Lata{\l}a's conjecture. This argument can be adapted
to the case of general parameter~$k$, and it yields that
Lata{\l}a's conjecture is equivalent to Theorem~2 of his paper~[6].

More generally, it can be said that even if we cannot prove
inequality~(5.12), the bound we can give for the high moments of
a random polynomial $Z$ defined in~(5.1) depends on what kind of
estimate we can prove for the expression at the left-hand side 
of~(5.12). But the estimation of such an expression is a hard 
problem. The analogous problem in formula~(5.11) was much simpler.

The estimation of the left-hand side of (5.12) is much harder
than that of formula~(5.11), because in this case the supremum
of random trilinear forms (and not of random bilinear forms as 
in formula~(5.11)) has to be considered. Lata{\l}a tried to get 
a good estimate for such an expression by means of a good bound 
on a quantity denoted by $N(X,\rho_\alpha,\e)$. The definition 
of this quantity was explained also in Section~2 of this paper. 
But his proof of the estimate for $N(\cdot)$ contained a serious 
error. I do not know whether the estimation of the quantity
$N(X,\rho_\alpha,\e)$ is the right way to give a good bound on
the expression in formula~(5.12). But the proof of such an 
estimate demands a deeper analysis than the method of paper~[6]. 
It should exploit the finer structure of the model we consider.

I do not know whether relation (5.12) is true. I can neither
prove it nor can I give a counter-example. I can only show with 
the help of the results in the present paper that a weaker form 
of the estimate~(5.12) holds with an upper bound $CM^{-1/4}$ 
instead of $CM^{-1/2}$. I briefly explain this.

Let us consider the random polynomial $Z$~in (5.1) with the 
difference that in the upper bounds of conditions (5.3) and (5.4) 
the numbers $R$ and $R^2$, $0\le R\le 1$, appear instead of 
$M^{-1/2}$ and $M^{-1}$ respectively. (In general, the number $R$ 
will take the same role in the next consideration as the parameter  
$M^{-1/2}$ did before.) The statement of relation (5.7), and as a 
consequence of relation~(5.8) can be reversed in the following way. 
They remain valid if the less or equal sign is replaced by the 
greater or equal sign in them, and the sufficiently large universal 
constant $C>0$ is replaced by another sufficiently small universal 
constant $C>0$. (See e.g]~[5]).

The random polynomial $Z$ defined in~(5.1) satisfies the 
inequality $EZ^{2\bar M}\le C_1^{\bar M}\bar M^{\bar M}$ with 
$\bar M=\frac1R$ and a sufficiently large constant $C_1>0$ if the 
modified version of (5.3) and (5.4) holds (with the replacement of
$M^{-1/2}$ by $R$). This follows from the results of these paper,
e.g. from formula~(2.11) with $k=3$. This estimate together with 
the above mentioned reversed form of formula~(5.8) or of its 
equivalent version given in~(5.10) with parameter $\bar M$ instead 
of the parameter~$M$ and the H\"older inequality imply that
$$
\align
C_1^{\bar M}\bar M^{\bar M}\ge EZ^{2\bar M}
\ge C_1^{\bar M}\bar M^{2\bar M} 
&\ge  C^{\bar M} \bar M^{2\bar M} E\(\sup\Sb v(j),\;w(k)\\ 
\summ v^2(j)\le1,\;\; \summ w^2(k)\le 1\endSb Y(v,w)\)^{2\bar M} \\
&\ge  C^{\bar M} \bar M^{2\bar M} \( E\left|\sup\Sb v(j),\;w(k)\\ 
\summ v^2(j)\le1,\;\; \summ w^2(k)\le 1\endSb Y(v,w)\right|\)^{2\bar M}.
\endalign
$$
The above estimates imply that
$$
E\left|\sup\Sb v(j),\;w(k)\\ 
\summ v^2(j)\le1,\;\; \summ w^2(k)\le 1\endSb Y(v,w)\right|\le C
\bar M^{-1/2}=CR^{1/2}
$$
with an appropriate constant $C>0$. Put $R=M^{-1/2}$, $\bar M=M^{1/2}$.
With such a choice we get the weakened form of relation (5.12) with the
bound $CM^{-1/4}$ instead of $CM^{-1/2}$ on its right-hand side.

I also remark that this estimate together with the results of this
section supply a slightly better bound on $EZ^{2M}$, than the bound
supplied by the results in Section~2. Namely 
$EZ^{2M}\le \(\max(M,M^2R,M^3R^4)\)^M$ instead of
$EZ^{2M}\le \(\max(M,M^3R^2)\)^M$ if the modified version of 
relations~(5.2), (5.3) and~(5.4) hold with the replacement of $M^{-1/2}$
by $R$ in their upper bound.

\beginsection

{\bf References:}

\medskip


\item{1.)} de la Pe\~na, V. H. and  Montgomery--Smith, S. (1995)
Decoupling inequalities for the tail-probabilities of multivariate
$U$-statistics. {\it Ann. Probab.}, {\bf 23}, 806--816

\item{2.)} Dobrushin, R. L. (1979) Gaussian and their subordinated
fields.  {\it Annals of Probability}\/ {\bf 7}, 1-28

\item{3.)} Hanson, D. L. and Wright, F. T. (1971) A bound on the
tail probabilities for quadratic forms  in independent random
variables. {\it Ann. Math. Statist.} {\bf 42} 52--61

\item{4.)} It\^o K. (1951) Multiple Wiener integral. {\it J. Math.
Soc. Japan}\/  {\bf3}.  157--164

\item{5.)} Lata\l{a}, R. (1999) Tail an moment estimates for some
types of chaos. {\it Studia Math.} {\bf 135} 39--53

\item{6.)} Lata\l{a}, R. (2006) Estimates of moments and tails of
Gaussian chaoses. {\it Annals of Probability} {\bf34} 2315--2331

\item{7.)} Ledoux, M. (2001) The concentration of measure phenomenon.
American Mathematical Society, Providence, RI

\item{8.)} Major, P. (1981) Multiple Wiener--It\^o integrals. {\it
Lecture Notes in Mathematics\/} {\bf 849}, Springer Verlag, Berlin,
Heidelberg, New York,

\item{9.)} Major, P. (2006) A multivariate generalization of
Hoeffding's inequality. {\it Electronic Communication in
Probability} {\bf 2} (220--229)

\item{10.)} Major, P. (2007) On a multivariate version of
Bernstein's inequality. {\it Electronic Journal of Probability}\/
{\bf12} 966--988

\item{11.)} Major, P. (2008) On the estimation of degenerate
$U$-statistics.  (in preparation.)

\item{12.)} Malyshev, V. A. and Minlos, R. A. (1991) Gibbs Random
Fields. Method of cluster expansion. Kluwer, Academic Publishers,
Dordrecht

\bye